\documentclass[ejsv2,preprint,noshowframe]{imsart}
\RequirePackage[numbers]{natbib}
\RequirePackage[colorlinks,citecolor=blue,urlcolor=blue]{hyperref}
\RequirePackage{graphicx}
\RequirePackage{bm}
\RequirePackage{multirow}
\usepackage{algorithm}
\usepackage{algpseudocode}
\usepackage{mathtools}
\usepackage{float}

\startlocaldefs

\theoremstyle{plain}

\newtheorem{theorem}{Theorem}[section]
\newtheorem{lemma}[theorem]{Lemma}

\theoremstyle{definition}

\theoremstyle{remark}
\newtheorem{remark}{Remark}[section]

\endlocaldefs

\begin{document}
\begin{frontmatter}
\title{Sharp minimax rates and adaptive estimation of latent proportions in bivariate Poisson mixtures}

\runtitle{Minimax estimation of latent proportions}

\begin{aug}
\author{}
\author[A,B]{\snm{Essoham Ali}\ead[label=e1]{essoham.ali@univ-ubs.fr}}
\address[A]{Institut de Math\'ematiques Appliqu\'ees, Universit\'e Catholique de l'Ouest, Angers, France.\printead[presep={\ }]{}}
 
\address[B]{Universit\'e Bretagne Sud, UMR CNRS 6205, LMBA, Vannes, France.\\ \printead[presep={\ }]{e1}}
\runauthor{E. Ali}
\end{aug}
\begin{abstract}
We study estimation of the mean latent composition in a
bivariate Poisson mixture. Conditional on an unobserved
intensity $S$ and proportion $U$, two counts are independent
Poisson variables with respective means $SU$ and $S(1-U)$.
The mixing law is unrestricted on $[a,b]\times[0,1]$,
where $0<a<b<\infty$, and dependence between $S$ and $U$
is allowed. Although the mixing distribution is identifiable,
the usual ratio of mean counts estimates an intensity-weighted
proportion and may be biased for $\theta=\mathbb E[U]$.
We construct polynomial estimators based on factorial moments
and establish the exact leading exponent of the minimax
squared-error risk on the logarithmic scale, without claiming
an exact multiplicative minimax constant.
We also derive the sharp logarithmic separation exponent
for testing zero latent covariance and construct an adaptive
estimator achieving parametric squared-error risk on the
zero-covariance submodel while retaining the optimal
full-model logarithmic minimax exponent.
Simulations illustrate the bias-variance trade-off,
the role of polynomial regularization, and the distinction
between theoretically guaranteed and empirically calibrated
adaptation.
\end{abstract}

\begin{keyword}[class=MSC]
\kwd[Primary ]{62G05}
\kwd[; secondary ]{62G20}
\kwd{62F07}
\end{keyword}

\begin{keyword}
\kwd{Poisson mixtures}
\kwd{Latent proportions}
\kwd{Statistical inverse problems}
\kwd{Sharp logarithmic minimax rates}
\kwd{Polynomial approximation}
\kwd{Adaptive estimation}
\kwd{Covariance testing}
\end{keyword}

\end{frontmatter}

\section{Introduction}\label{sec:introduction}

Count data collected from heterogeneous observational units
often reflect two sources of variation: the overall intensity
of event occurrence and the relative composition of different
event types. Examples arise in ecological sampling, recurrent
event studies, and biological count data. When both intensity
and composition vary across units, the average composition
of the units need not coincide with the composition obtained
by pooling their observed counts. This distinction is
particularly important when latent intensity and composition
are statistically dependent.

Mixed Poisson models provide a flexible framework for
representing unobserved heterogeneity in count data
\cite{KarlisXekalaki2005}. Multivariate extensions have
also been investigated to describe dependence between
count components \cite{AitchisonHo1989,Holgate1964}.
In this paper, we consider a bivariate Poisson splitting
model in which the total latent intensity and its allocation
between two count components may be arbitrarily dependent.
Our objective is to characterize the fundamental limits
of estimating the mean latent composition and to investigate
whether estimation can adapt to a simpler submodel defined
by a zero-covariance restriction.

Let $Y_i=(X_{i1},X_{i2})$, $i=1,\ldots,n$, be independent
and identically distributed observations. Conditionally on
latent variables $(S_i,U_i)=(s,u)$, suppose that
\begin{equation}\label{eq:intro-model}
X_{i1}\mid(s,u)\sim\operatorname{Pois}(su),
\qquad
X_{i2}\mid(s,u)\sim
\operatorname{Pois}\{s(1-u)\},
\end{equation}
independently. Here $S_i$ represents the total latent
intensity, while $U_i$ denotes the proportion allocated
to the first count component. The joint mixing distribution
$Q$ is unrestricted on $[a,b]\times[0,1]$, where
$0<a<b<\infty$ are fixed and known. In particular,
we impose neither independence nor smoothness assumptions
on the latent variables.

The parameter of interest is the mean latent composition
$\theta(Q)=\mathbb E_Q[U]$, which describes the average
composition of a randomly selected observational unit.
This parameter differs from the proportion of first-type
events in the pooled population. Indeed, writing
$D=X_1+X_2$, we have
\begin{equation}\label{eq:intro-ratio}
\frac{\mathbb E_Q[X_1]}{\mathbb E_Q[D]}
=
\theta(Q)+
\frac{\operatorname{Cov}_Q(S,U)}
{\mathbb E_Q[S]}.
\end{equation}
Consequently, the conventional ratio of aggregate counts
may be inconsistent for the mean latent composition.
The discrepancy is governed exactly by the covariance
between total intensity and composition. Conversely,
when this covariance vanishes, the ratio consistently
estimates $\theta(Q)$, even if the latent variables
remain dependent. This observation motivates both
estimation over the unrestricted mixing class and
inference on the covariance-null submodel.

Although the joint mixing distribution is identifiable
under the stated support assumptions, identifiability
alone does not determine the accuracy with which
$\theta(Q)$ can be estimated. Recovering the unweighted
mean composition requires separating the effect of
latent intensity from the allocation of counts between
the two components. This leads to an inverse problem
whose statistical difficulty differs substantially
from that of estimating the observable intensity-weighted
composition.

The estimation of mixing distributions and their
functionals has a long history in mathematical statistics
\cite{KieferWolfowitz1956,Lindsay1983,Lindsay1995}.
More generally, optimal recovery and minimax theory
have shown that the estimation of a functional can
exhibit convergence rates substantially different
from those associated with recovery of the underlying
distribution \cite{Donoho1994,CaiLow2004}.

A closely related contribution is that of
Cavalier and Hengartner
\cite{CavalierHengartner2009}, who study
the estimation of linear functionals of
mixing distributions from univariate
Poisson mixture observations.
Their analysis develops inversion methods
for analytic functionals, Gaussian
approximations for less regular
functionals, and estimation procedures
under additional regularity assumptions.
These results provide an important
methodological foundation for functional
estimation in Poisson inverse problems.

The present problem differs in both
its observation structure and its target
functional. Although the total count
$D=X_1+X_2$ follows a univariate
Poisson mixture, its distribution depends
only on the marginal law of $S$ and
does not identify
$\theta(Q)=\mathbb E_Q[U]$
under arbitrary dependence between
$S$ and $U$.
The joint observations instead provide
mixed factorial moments satisfying
\[
\mathbb E_Q\left[
X_1(D-1)_j\mathbf 1_{\{D\ge1\}}
\right]
=
\mathbb E_Q[S^{j+1}U],
\quad j\ge0.
\]
Consequently, estimation of the mean
latent composition can be formulated
through polynomial approximation of
the reciprocal intensity $1/s$
on $[a,b]$.

Equivalently, introducing
$\Lambda_1=SU$ and
$\Lambda_2=S(1-U)$
represents the model as a bivariate
Poisson mixture, with target functional
$\mathbb E_Q[
\Lambda_1/(\Lambda_1+\Lambda_2)]$.
Thus, the problem remains within
the general framework of functional
estimation in mixture models.
Its specific difficulty arises from
the joint observation operator and
the reciprocal dependence on total
latent intensity.

Our contribution is to establish
matching minimax upper and lower bounds
at the level of the leading logarithmic
exponent for this bivariate experiment
over the unrestricted class of mixing
distributions supported on
$[a,b]\times[0,1]$.
In particular, we identify the exact
leading logarithmic minimax exponent,
including its explicit dependence
on the support endpoints.
The lower bound controls the total
variation distance between the full
bivariate observation laws and
therefore does not follow merely
by discarding component labels
and applying a univariate result.
We further investigate the corresponding
covariance-testing problem and adaptation
to the covariance-null submodel,
neither of which is addressed in
\cite{CavalierHengartner2009}.

Adaptive estimation in Poisson mixture models
has also been investigated by \cite{Hengartner1997},
particularly for mixing-density estimation over
smoothness classes. Our adaptation problem has
a different structure: the relevant classes are
distinguished by a restriction on latent covariance
rather than by smoothness. The objective is to
exploit the parametric behavior of the ratio
estimator when this covariance vanishes, without
sacrificing the minimax logarithmic exponent
over the unrestricted class. This question is
also connected to the broader literature on
adaptive estimation and estimator comparison
\cite{Lepski1991,GoldenshlugerLepski2011}.

Recent developments in nonparametric Poisson mixture
models include the estimation of mixing distributions
under smoothed optimal transport losses
\cite{LimHan2024} and minimum-distance approaches
to empirical Bayes estimation
\cite{OptimalEB2025,FavaroFortiniJana2026}.
These studies address inferential objectives and
loss functions that differ from those considered
here. In contrast, our analysis focuses on the
minimax estimation of mean latent composition
in bivariate Poisson mixtures, the detection
of latent covariance, and adaptation to a
covariance-null submodel.

The first main contribution of this paper is
a sharp characterization of the minimax difficulty
of estimating the mean latent composition.
Using a polynomial estimator based on factorial
moments and a matching lower-bound construction
derived from Chebyshev approximation, we establish
\begin{equation}\label{eq:intro-minimax}
\lim_{n\to\infty}
\frac{-\log R_n^*(a,b)}
{\log n/\log\log n}
=
2\log\left(
\frac{\sqrt b+\sqrt a}
{\sqrt b-\sqrt a}
\right),
\end{equation}
where $R_n^*(a,b)$ denotes the minimax
squared-error risk over the unrestricted
mixing class. The result determines the exact
leading coefficient on the logarithmic scale
and makes explicit its dependence on the
support bounds of the latent intensity.
The sharpness concerns the logarithmic exponent,
not an exact multiplicative minimax constant.

The second main contribution concerns the
statistical detectability of latent covariance.
We establish a sharp logarithmic separation
exponent for testing the covariance-null
restriction against unrestricted alternatives.
The proof combines uniform covariance estimation
with a construction of mixing distributions
that have separated latent covariances but
induce asymptotically indistinguishable
count distributions. 
This result quantifies the intrinsic
difficulty of detecting the latent
covariance responsible for the bias
of the aggregate-count ratio. 

The third main contribution is an adaptive
estimator combining the conventional ratio
with a robust polynomial estimator.
We show that it achieves parametric
squared-error risk of order $n^{-1}$
uniformly over the covariance-null submodel,
while retaining the optimal logarithmic
minimax exponent over the unrestricted class.
The construction relies on estimator comparison
and uniform concentration inequalities,
rather than on consistent identification
of the covariance-null submodel.
The resulting adaptation guarantee concerns
risk and does not assert semiparametric
efficiency or uniformly accurate
finite-sample model selection.

As a preliminary theoretical result,
we establish identifiability of the
unrestricted joint mixing distribution
from the observed bivariate counts.
The strictly positive lower bound on
latent intensity is essential to the
present formulation, while compact
support determines the polynomial
approximation problem underlying
the minimax analysis.

Numerical experiments investigate the
finite-sample behavior of the proposed
estimators under different forms of
latent dependence. We examine the effects
of polynomial degree, support specification,
adaptive threshold calibration, and covariance
testing under fixed and local alternatives.
The experiments also highlight the
conservativeness of theoretically calibrated
selection thresholds at moderate sample sizes.
Accordingly, we distinguish the procedure
covered by the uniform theoretical guarantees
from an empirically calibrated implementation.

The remainder of the paper is organized as
follows. Section~\ref{sec:model} presents
the model and identifiability result.
Section~\ref{sec:minimax} develops polynomial
estimation and sharp minimax bounds.
Section~\ref{sec:adaptation} investigates
covariance testing and adaptive estimation.
Section~\ref{sec:sim} reports numerical
experiments. Section~\ref{sec:discussion}
discusses the main findings, limitations,
and directions for future research.
Proofs and auxiliary results are provided
in the appendix. The R code and numerical results used in this study
are publicly available at
\url{https://github.com/essohali/EJS-Latent-Proportions}.

\section{Model and identifiability}\label{sec:model}

\subsection{Latent Poisson splitting model and observable moments}
\label{subsec:model}

Let $Y_i=(X_{i1},X_{i2})$, $i=1,\ldots,n$, be
independent and identically distributed random vectors
taking values in $\mathbb N_0^2$. For each $i$, let
$(S_i,U_i)$ be an unobserved random vector with common
distribution $Q$ supported on $\mathcal K_{a,b}=[a,b]\times[0,1],
\quad 0<a<b<\infty.$
Conditional on $(S_i,U_i)=(s,u)$, the component counts
are independent and satisfy
\begin{equation}\label{eq:conditional-poisson}
X_{i1}\mid(s,u)\sim\operatorname{Pois}(su),
\quad
X_{i2}\mid(s,u)\sim\operatorname{Pois}\{s(1-u)\}.
\end{equation}
No independence assumption is imposed between $S_i$
and $U_i$. We denote by $\mathcal Q(a,b)$ the set of
all Borel probability measures on $\mathcal K_{a,b}$.
Our target is the mean latent composition
\begin{equation}\label{eq:theta}
\theta(Q)
=
\int_{\mathcal K_{a,b}}u\,Q(ds,du)
=
\mathbb E_Q[U].
\end{equation}

The conditional construction in
\eqref{eq:conditional-poisson} is an instance of
Poisson splitting \cite{Kingman1993}; see also
\cite[Chapter~VI]{Feller1968}. Classical work on
the identifiability of mixture distributions includes
\cite{Teicher1960,Teicher1961}, while the estimation
of linear functionals of univariate Poisson mixing
distributions is investigated by
\cite{CavalierHengartner2009}. These works provide
the general statistical context for the model.
The identifiability result established below is
used as a preliminary step toward the minimax
analysis rather than as a separate contribution
to the general theory of mixture identifiability.

For a generic observation $Y=(X_1,X_2)$, the joint
probability mass function is
\begin{equation}\label{eq:mixture-pmf}
p_Q(k,\ell)
=
\frac{1}{k!\ell!}
\int_{\mathcal K_{a,b}}
e^{-s}s^{k+\ell}u^k(1-u)^\ell
\,Q(ds,du),
\quad (k,\ell)\in\mathbb N_0^2.
\end{equation}
Writing $D=X_1+X_2$, the Poisson--binomial
splitting identity gives
\begin{equation}\label{eq:splitting}
D\mid(S=s,U=u)\sim\operatorname{Pois}(s),
\quad
X_1\mid(D=d,S=s,U=u)\sim\operatorname{Bin}(d,u),
\end{equation}
for every $d\in\mathbb N_0$. Thus the distribution
of $D$ depends only on the marginal law of $S$,
whereas the joint count distribution also contains
information about the latent composition.

The first observable moments satisfy
\begin{equation}\label{eq:first-moments}
\mathbb E_Q[D]
=
\mathbb E_Q[S]=:\mu_Q,
\quad
\mathbb E_Q[X_1]
=
\mathbb E_Q[SU]=:\eta_Q.
\end{equation}
Consequently,
\begin{equation}\label{eq:ratio-bias}
\frac{\mathbb E_Q[X_1]}{\mathbb E_Q[D]}
=
\theta(Q)
+
\frac{\operatorname{Cov}_Q(S,U)}{\mu_Q},
\quad a\leq\mu_Q\leq b.
\end{equation}
The ratio of observable mean counts therefore
represents the intensity-weighted composition
$\eta_Q/\mu_Q$, which generally differs from the
unweighted target $\theta(Q)$. Equality holds
whenever $\operatorname{Cov}_Q(S,U)=0$.
Independence of $S$ and $U$ is sufficient but
not necessary for this condition. This distinction
will play a central role in the covariance testing
and adaptive estimation problems considered
in Section~\ref{sec:adaptation}.

For integers $j\geq0$, define the falling factorial
by $(x)_0=1$ and
$(x)_j=x(x-1)\cdots(x-j+1)$ for $j\geq1$.
Conditional Poisson factorial-moment identities give
\begin{equation}\label{eq:factorial-identities}
\begin{aligned}
\mathbb E_Q[(D)_j]
&=
\mathbb E_Q[S^j],\\
\mathbb E_Q\!\left[
X_1(D-1)_j\mathbf 1_{\{D\geq1\}}
\right]
&=
\mathbb E_Q[S^{j+1}U].
\end{aligned}
\end{equation}
The expression
$X_1(D-1)_j\mathbf 1_{\{D\geq1\}}$
is defined to be zero when $D=0$.
Indeed, conditionally on $(D,S,U)$,
\[
\mathbb E[X_1\mid D,S,U]=DU.
\]
For $D\geq1$, the identity
$D(D-1)_j=(D)_{j+1}$ and the conditional
factorial moments of the Poisson distribution
yield the second equality in
\eqref{eq:factorial-identities}.

These identities provide the link between
observable count statistics and latent moments.
In particular, they make it possible to estimate
$\mathbb E_Q[U]$ through polynomial approximation
of the reciprocal function on $[a,b]$, as developed
in Section~\ref{sec:minimax}.

\subsection{Identifiability of the mixing distribution}
\label{subsec:identifiability}

Define the latent intensity vector $(\Lambda_1,\Lambda_2)
=
(SU,S(1-U)).$
The mapping $\Phi(s,u)=(su,s(1-u))$
is a continuous bijection from
$\mathcal K_{a,b}$ onto the compact set
\[
\mathcal L_{a,b}
=
\left\{
(\lambda_1,\lambda_2)\in[0,\infty)^2:
a\leq\lambda_1+\lambda_2\leq b
\right\}.
\]
Its inverse is continuous and is given by $s=\lambda_1+\lambda_2$
\text{and}
$u=\frac{\lambda_1}{\lambda_1+\lambda_2}.
$
The condition $a>0$ ensures that the denominator
is uniformly bounded away from zero.

The identifiability argument relies on two
properties of the model. First, the joint
factorial moments of the observed counts identify
the mixed moments of the latent intensity vector:
\[
\mathbb E_Q[(X_1)_k(X_2)_\ell]
=
\mathbb E_Q[\Lambda_1^k\Lambda_2^\ell],
\quad k,\ell\geq0.
\]
Second, since $\mathcal L_{a,b}$ is compact,
polynomials in $(\lambda_1,\lambda_2)$ are dense
in the space of continuous functions on this
set. Consequently, the mixed moments uniquely
determine the distribution of
$(\Lambda_1,\Lambda_2)$, and the invertibility
of $\Phi$ then identifies $Q$.

\begin{theorem} 
\label{thm:identifiability}
Fix $0<a<b<\infty$. If
$Q,Q'\in\mathcal Q(a,b)$ induce the same
distribution of $(X_1,X_2)$ under
\eqref{eq:conditional-poisson}, then $Q=Q'$.
Consequently, $\theta(Q)$ is identifiable
from the joint count distribution.
\end{theorem}

\noindent\emph{Proof.}
See Appendix~\ref{app:proofs}.

\begin{remark} 
\label{rem:zero}
The strict lower bound $a>0$ excludes an
unobservable composition at zero total intensity.
More generally, if a model allows
$\Pr(S=0)>0$ without restricting the conditional
law of $U$ on $\{S=0\}$, then
$(X_1,X_2)=(0,0)$ almost surely on that stratum,
regardless of $U$. Changing the conditional law
of $U$ given $S=0$ can therefore change
$\mathbb E[U]$ without changing the observable
distribution. Thus identification of the
unrestricted latent composition law, and even
of its mean, can fail.

This does not rule out identification in
restricted submodels allowing $S=0$.
The upper bound $b<\infty$ ensures moment
determinacy through compactness in the
argument above; no necessity claim is made
for bounded support in general Poisson mixtures.
\end{remark}

\begin{remark}  
\label{rem:total}
Under the compact-support assumption, the
distribution of $D$ identifies the marginal
law of $S$: its factorial moments determine
all moments of $S$, and polynomials are dense
in $C([a,b])$. However, the distribution of
$D$ alone does not identify $\theta(Q)$
over $\mathcal Q(a,b)$.

For example, fixing any distribution of $S$
and taking either $U\equiv0$ or $U\equiv1$
leaves the distribution of $D$ unchanged,
while changing $\theta(Q)$ from zero to one.
Theorem~\ref{thm:identifiability} instead
uses the full joint count distribution.

Identifiability is a population-level property
and does not, by itself, imply a particular
convergence rate or numerical stability
for estimating $\theta(Q)$.
\end{remark}

Theorem~\ref{thm:identifiability} establishes
that the target functional is uniquely determined
by the observable joint distribution. This
result is preparatory: the principal statistical
question is not whether $\theta(Q)$ is
identifiable, but how accurately it can be
estimated from a finite sample.

In particular, the inverse relationship between
observable factorial moments and the latent
functional may be statistically ill-conditioned,
even though the mixing distribution is uniquely
identified. Section~\ref{sec:minimax} quantifies
this difficulty by establishing sharp logarithmic
minimax bounds for estimating $\theta(Q)$ over
the unrestricted mixing class.

\section{Sharp logarithmic minimax bounds}
\label{sec:minimax}

\subsection{Risk and polynomial construction}
\label{subsec:minimax-upper}

For the model introduced in Section~\ref{sec:model},
let $P_Q$ denote the distribution of a single
observation under $Q\in\mathcal Q(a,b)$.
We consider the minimax squared-error risk
\begin{equation}\label{eq:minimax-risk}
R_n^*(a,b)
=
\inf_{\widehat\theta}
\sup_{Q\in\mathcal Q(a,b)}
\mathbb E_Q
\left[
(\widehat\theta-\theta(Q))^2
\right],
\end{equation}
where the infimum is taken over all measurable
estimators based on $n$ independent observations.

The estimation of linear functionals of Poisson
mixing distributions is closely related to
statistical inverse problems; see
\cite{CavalierHengartner2009}.
Our analysis combines polynomial approximation
with the classical two-point method for minimax
lower bounds \cite{Tsybakov2009}.
The extremal properties of Chebyshev polynomials
play a central role in both constructions;
see \cite{Cheney1966}.

Define
\[
z=\frac{a+b}{b-a}>1,
\quad
r=\frac{\sqrt b-\sqrt a}{\sqrt b+\sqrt a}
\in(0,1)
\quad\text{and}\quad
\lambda=-\log r.
\]
Let $T_m$ denote the Chebyshev polynomial of
the first kind, characterized by
$T_m(\cos t)=\cos(mt)$.
For each integer $m\geq1$, define
\[
P_m(s)
=
\frac{
T_m\left(\frac{a+b-2s}{b-a}\right)
}{
T_m(z)
},
\quad
q_{m-1}(s)
=
\frac{1-P_m(s)}{s}
=
\sum_{j=0}^{m-1}c_{j,m}s^j.
\]
Since $P_m(0)=1$, the function $q_{m-1}$
is a polynomial of degree at most $m-1$.

The Chebyshev extremal property gives
\[
\sup_{s\in[a,b]}
|1-sq_{m-1}(s)|
=
\frac{1}{T_m(z)}
=
\frac{2r^m}{1+r^{2m}}
\leq 2r^m.
\]
Thus $q_{m-1}$ approximates the reciprocal
function $1/s$ with geometrically decreasing
weighted approximation error on $[a,b]$.

For $D_i=X_{i1}+X_{i2}$, define
\[
H_j(Y_i)
=
\begin{cases}
X_{i1}(D_i-1)_j, & D_i\geq1,\\
0, & D_i=0,
\end{cases}
\quad j\geq0,
\]
where $(x)_j$ denotes the falling factorial.
In particular, $H_0(Y_i)=X_{i1}$.
We consider the polynomial estimator
\begin{equation}\label{eq:poly-estimator}
Z_{i,m}
=
\sum_{j=0}^{m-1}c_{j,m}H_j(Y_i),
\qquad
\widehat\theta_{n,m}
=
\frac1n\sum_{i=1}^n Z_{i,m}.
\end{equation}

By the factorial-moment identities in
\eqref{eq:factorial-identities},
\[
\mathbb E_Q[Z_{i,m}]
=
\mathbb E_Q[SUq_{m-1}(S)].
\]
Consequently,
\[
\begin{aligned}
\left|
\mathbb E_Q[\widehat\theta_{n,m}]
-\theta(Q)
\right|
&=
\left|
\mathbb E_Q[
U\{Sq_{m-1}(S)-1\}
]
\right|\\
&\leq
\sup_{s\in[a,b]}
|sq_{m-1}(s)-1|
\leq 2r^m.
\end{aligned}
\]
The approximation bias therefore decreases
geometrically with the polynomial degree,
uniformly over $\mathcal Q(a,b)$.

Increasing the polynomial degree, however,
also amplifies sampling variability.
The coefficients of $q_{m-1}$ and the
factorial-moment statistics entering
$Z_{i,m}$ may grow rapidly with $m$.
A uniform second-moment bound, established
in Appendix~\ref{app:proofs}
(Lemma~\ref{lem:minimax-upper}), yields
\[
\sup_{Q\in\mathcal Q(a,b)}
\mathbb E_Q
\left[
(\widehat\theta_{n,m}-\theta(Q))^2
\right]
\leq
4r^{2m}
+
\frac{C_0C_1^m m!}{n},
\]
where $C_0,C_1>0$ depend only on $(a,b)$.

This inequality makes explicit the
bias--variance trade-off underlying
polynomial regularization.
The first term decreases geometrically,
whereas the second reflects the factorial
growth associated with estimating
higher-order latent moments.
In particular, taking
\[
m_n=
\left\lfloor
\frac{\log n}{\log\log n}
\right\rfloor
\]
makes the variance term negligible relative
to $r^{2m_n}$ on the logarithmic scale.
The resulting upper bound determines one
direction of the sharp minimax exponent.

\subsection{Minimax lower bound and sharp exponent}
\label{subsec:minimax-lower}

We next establish that the logarithmic
exponent obtained from the polynomial
construction cannot be improved uniformly
over $\mathcal Q(a,b)$.
The argument combines the duality of
uniform polynomial approximation with
a two-point testing construction.

Let $V_m
=
\operatorname{span}
\{s,s^2,\ldots,s^m\}
\subset C([a,b]).$
The Chebyshev extremal property implies
\[
\inf_{p\in V_m}
\sup_{s\in[a,b]}|1-p(s)|
=
\frac{1}{T_m(z)}.
\]
By duality, this approximation error can
be represented by a signed measure
$\nu_m$ on $[a,b]$ satisfying
\[
\|\nu_m\|_{\mathrm{var}}=1,
\qquad
\int s^j\,\nu_m(ds)=0,
\quad j=1,\ldots,m,
\]
and
\[
\int\nu_m(ds)
=
\frac{1}{T_m(z)}.
\]

Writing $\nu_m=\nu_m^+-\nu_m^-$
for its Jordan decomposition, define
two mixing distributions by
\[
\begin{aligned}
Q_m^+(ds,du)
&=
\nu_m^+(ds)\delta_1(du)
+
\nu_m^-(ds)\delta_0(du),\\
Q_m^-(ds,du)
&=
\nu_m^+(ds)\delta_0(du)
+
\nu_m^-(ds)\delta_1(du).
\end{aligned}
\]
Both distributions belong to
$\mathcal Q(a,b)$ and have the same
marginal distribution of $S$.
Nevertheless, their target functionals
satisfy
\[
\theta(Q_m^+)-\theta(Q_m^-)
=
\frac{1}{T_m(z)}.
\]

Let $P_m^+$ and $P_m^-$ denote the
corresponding distributions of one
observed bivariate count vector.
Because $U$ takes values in $\{0,1\}$
under these mixing distributions,
the observations are supported on
the coordinate axes.
The moment cancellation conditions,
combined with a Taylor expansion and
a Poisson tail bound, yield
\[
d_{\mathrm{TV}}(P_m^+,P_m^-)
\leq
\frac{(2b)^{m+1}}{(m+1)!}.
\]
Importantly, this bound concerns the
full bivariate observation distributions,
not merely the distributions of the
total counts.

Applying the two-point minimax
inequality and the product total-variation
bound gives
\[
R_n^*(a,b)
\geq
\frac{1}{8T_m(z)^2}
\left[
1-
\frac{n(2b)^{m+1}}{(m+1)!}
\right]_+,
\]
where $[x]_+=\max\{x,0\}$.
The complete construction and proof
are given in Appendix~\ref{app:proofs}
(Lemma~\ref{lem:minimax-lower}).

The lower bound complements the
polynomial risk upper bound.
Taking
\[
m_n=
\left\lceil
(1+\varepsilon)
\frac{\log n}{\log\log n}
\right\rceil,
\qquad \varepsilon>0,
\]
makes the total-variation contribution
asymptotically negligible.
Together with the upper-bound degree
choice, this identifies the exact
leading coefficient of the negative
logarithmic minimax risk.

\begin{theorem} 
\label{thm:minimax-rate}
For every fixed $0<a<b<\infty$,
\begin{equation*}
\lim_{n\to\infty}
\frac{-\log R_n^*(a,b)}
{\log n/\log\log n}
=
2\log\left(
\frac{\sqrt b+\sqrt a}
{\sqrt b-\sqrt a}
\right).
\end{equation*}
\end{theorem}

\noindent\emph{Proof.}
See Appendix~\ref{app:proofs}.

Theorem~\ref{thm:minimax-rate} characterizes
the intrinsic difficulty of estimating
the unweighted latent composition under
arbitrary dependence between intensity
and composition.
Although the proof uses classical tools
from polynomial approximation and minimax
testing, the resulting exponent is obtained
for the specific functional
$\theta(Q)=\mathbb E_Q[U]$ under the
bivariate Poisson splitting observation
model.

The distinction from univariate Poisson
mixture functional estimation is structural.
The total count $D$ identifies the marginal
law of $S$ but does not identify
$\theta(Q)$ over the unrestricted class.
The joint count observations instead
provide the mixed factorial moments
needed to estimate the target.
Moreover, the lower-bound construction
controls the total-variation distance
between the complete bivariate count
distributions.
The matching bounds therefore establish
the minimax exponent for the full
statistical experiment rather than
for a reduction based only on total counts.

\begin{remark}
\label{rem:minimax-sharpness}
Theorem~\ref{thm:minimax-rate} identifies
the exact coefficient of
$\log n/\log\log n$ in the negative
logarithm of the minimax squared-error risk.
It does not establish an asymptotic
equivalence of the form
\[
R_n^*(a,b)
\sim
C\exp(-2\lambda L_n),
\qquad
L_n=\frac{\log n}{\log\log n},
\]
nor does it determine an optimal
finite-sample multiplicative constant.
The support bounds $a$ and $b$ are
assumed known and fixed throughout
this section.
\end{remark}

\section{Detection of latent dependence and adaptive estimation}
\label{sec:adaptation}

\subsection{Testing zero latent covariance}
\label{subsec:testing}

The decomposition in \eqref{eq:ratio-bias} suggests
studying the covariance-null submodel $\mathcal Q_0(a,b)
=
\{Q\in\mathcal Q(a,b):C(Q)=0\},$
where $C(Q)
=
\operatorname{Cov}_Q(S,U)
=
\eta_Q-\mu_Q\theta(Q).$
The restriction $C(Q)=0$ is weaker than
independence between $S$ and $U$.
In particular, it allows nonlinear forms
of latent dependence.

Writing $\theta_R(Q)=\frac{\eta_Q}{\mu_Q},$
we have the identity
\begin{equation}\label{eq:covariance-functional}
C(Q)
=
\mu_Q\{\theta_R(Q)-\theta(Q)\}.
\end{equation}
Since $a\leq\mu_Q\leq b$, testing
zero covariance is equivalent to testing
whether the observable intensity-weighted
composition coincides with the mean
latent composition.

This identity also reveals a connection
between covariance detection and functional
estimation. The quantities $\eta_Q$ and
$\mu_Q$ are ordinary observable moments,
whereas $\theta(Q)$ is the functional
whose estimation requires inversion of
the Poisson mixture operator.
Consequently, the difficulty of detecting
small values of $C(Q)$ is driven by the
same latent functional that determines
the minimax estimation problem in
Section~\ref{sec:minimax}.

For a sequence $\rho_n>0$, consider
\begin{equation}\label{eq:covariance-test}
H_0:Q\in\mathcal Q_0(a,b)
\quad\text{against}\quad
H_1:Q\in\mathcal Q(a,b),\
|C(Q)|\geq\rho_n.
\end{equation}
For a measurable test $\phi_n\in\{0,1\}$,
define its maximal sum of errors by
\[
\begin{aligned}
\mathcal E_n(\phi_n,\rho_n)
={}&
\sup_{Q\in\mathcal Q_0(a,b)}
P_Q^{\otimes n}(\phi_n=1)\\
&+
\sup_{\substack{Q\in\mathcal Q(a,b)\\
|C(Q)|\geq\rho_n}}
P_Q^{\otimes n}(\phi_n=0),
\end{aligned}
\]
where a supremum over an empty alternative
is taken to be zero.
The distinction between estimation and
testing rates is standard in nonparametric
inference; see
\cite{IngsterSuslina2003,Tsybakov2009}.

Let
\[
\widehat\mu_n=\frac1n\sum_{i=1}^nD_i,
\qquad
\widehat\eta_n=\frac1n\sum_{i=1}^nX_{i1}.
\]
Let $\widehat\theta_{n,m}$ be the polynomial
estimator defined in \eqref{eq:poly-estimator},
and set $\overline\theta_{n,m}
=
\min\{1,\max\{0,\widehat\theta_{n,m}\}\}.$
We estimate the latent covariance by
\begin{equation}\label{eq:covariance-estimator}
\widehat C_{n,m}
=
\widehat\eta_n
-
\widehat\mu_n\overline\theta_{n,m}.
\end{equation}

The projection ensures that
$\overline\theta_{n,m}\in[0,1]$ and
does not increase its squared error
relative to $\theta(Q)$.
Using \eqref{eq:covariance-functional},
the boundedness of $\mu_Q$, and the
uniform polynomial risk bound, we obtain
\[
\sup_{Q\in\mathcal Q(a,b)}
\mathbb E_Q
\left[
\{\widehat C_{n,m}-C(Q)\}^2
\right]
\leq
K_0\left(
r^{2m}
+
\frac{K_1^m m!}{n}
+
\frac1n
\right),
\]
where $K_0,K_1>0$ depend only on $(a,b)$.
The proof is given in
Appendix~\ref{app:proofs},
Lemma~\ref{lem:covariance-estimation}.

This bound yields a test based on
thresholding $|\widehat C_{n,m}|$.
For polynomial degrees of order
$\log n/\log\log n$, the estimation
error is controlled on the logarithmic
scale by the approximation term.
A threshold slightly larger than
$r^{m}$ therefore permits uniform
testing against sufficiently separated
covariance alternatives.

The converse argument uses a pair
of mixing distributions that are
difficult to distinguish from their
observed bivariate counts.
Starting from the extremal signed
measure used in the minimax lower-bound
construction, one obtains two laws
$Q_m^+$ and $Q_m^-$ with the same
marginal distribution of $S$.
Their average $Q_{0,m}
=
\frac12(Q_m^++Q_m^-)$
belongs to the covariance-null submodel.
Indeed, under $Q_{0,m}$, the latent
composition is Bernoulli with parameter
$1/2$ and is independent of $S$.

By contrast, the moment cancellation
conditions imply
\[
C(Q_m^+)
=
-\frac{\mu_m E_m}{2},
\]
where
\[
\mu_m=\mathbb E_{Q_m^+}[S],
\qquad
E_m=\frac{1}{T_m(z)}.
\]
In particular,
\[
|C(Q_m^+)|
\geq
\frac{a}{2T_m(z)}.
\]
The same construction controls the
total-variation distance between the
corresponding observed bivariate
distributions.

Thus the testing lower bound is based
on alternatives with nonzero latent
covariance whose observable count
distributions remain close to a
covariance-null distribution.
This connects the detection boundary
to the Chebyshev approximation error
appearing in the minimax estimation
problem.

\begin{theorem} 
\label{thm:testing}
Let $L_n=\frac{\log n}{\log\log n}$ and
$\lambda=-\log r.$
For every fixed $\varepsilon\in(0,1)$,
there exist tests $\phi_n$ such that $\mathcal E_n\left(
\phi_n,
\exp\{-(1-\varepsilon)\lambda L_n\}
\right)
\longrightarrow 0.$
Conversely, for every fixed $\varepsilon>0$, $\inf_{\phi_n}
\mathcal E_n\left(
\phi_n,
\exp\{-(1+\varepsilon)\lambda L_n\}
\right)
\longrightarrow 1.$
Thus, $\lambda$ is the sharp logarithmic
separation exponent for testing zero
latent covariance over the unrestricted
mixing class.
\end{theorem}
The sharpness concerns the leading logarithmic
separation exponent; no exact separation constant
or boundary-case characterization is asserted.

\noindent\emph{Proof.}
See Appendix~\ref{app:proofs}.

Theorem~\ref{thm:testing} identifies
the separation scale at which latent
covariance becomes statistically
detectable from the observed counts.
The exponent $\lambda$ is one-half
of the exponent $2\lambda$ governing
the minimax squared-error risk in
Theorem~\ref{thm:minimax-rate}.
This relationship reflects the
difference between a separation
measured in the scale of the functional
and an estimation risk measured
in squared error.

The result concerns detection of
nonzero covariance, not detection
of arbitrary dependence between
$S$ and $U$. In particular,
distributions with dependent latent
variables may belong to
$\mathcal Q_0(a,b)$.

\subsection{Adaptation to the zero-covariance submodel}
\label{subsec:adaptive}

On $\mathcal Q_0(a,b)$, the target
equals the observable ratio
$\eta_Q/\mu_Q$. This motivates the
conventional ratio estimator
\[
\widehat\theta_{R,n}
=
\begin{cases}
\widehat\eta_n/\widehat\mu_n,
&\widehat\mu_n>0,\\
0,&\widehat\mu_n=0.
\end{cases}
\]
Since $0\leq X_{i1}\leq D_i$,
this estimator takes values in $[0,1]$.

The ratio estimator has uniformly
parametric squared-error risk for
estimating $\theta_R(Q)$:
\[
\sup_{Q\in\mathcal Q(a,b)}
\mathbb E_Q
\left[
\{\widehat\theta_{R,n}
-\theta_R(Q)\}^2
\right]
\leq\frac{K_R}{n},
\]
for a constant $K_R>0$ depending
only on $(a,b)$; see
Appendix~\ref{app:proofs},
Lemma~\ref{lem:ratio-risk}.
Consequently, it achieves risk
of order $n^{-1}$ for estimating
$\theta(Q)$ on $\mathcal Q_0(a,b)$.

Outside this submodel, however,
the ratio estimator generally
converges to $\theta_R(Q)$ rather
than $\theta(Q)$. Its asymptotic
bias is therefore determined by
$C(Q)/\mu_Q$.

The adaptive estimation problem is
to exploit the parametric behavior
of the ratio estimator when
$C(Q)=0$, while preserving the
optimal logarithmic minimax
exponent over $\mathcal Q(a,b)$.

To this end, we combine the ratio
estimator with a robust polynomial
estimator. Median-of-means methods
provide exponential deviation
control under finite second-moment
assumptions; see
\cite{DevroyeLerasleLugosiOliveira2016}.

Set $m_n=
\max\{1,\lfloor L_n\rfloor\},
\quad
K_n=\lceil c\log n\rceil,$
where $c>0$ is a sufficiently large
fixed constant. For sufficiently
large $n$, partition the sample
into $K_n$ nonempty blocks of sizes
differing by at most one.

For each block, compute the mean
of $Z_{i,m_n}$ defined in
\eqref{eq:poly-estimator}.
Let $\widetilde\theta_{P,n}$
be the projection onto $[0,1]$
of the median of these block means.

Let $V_m=V_0V_1^m m!$
be the uniform variance bound
established in Appendix~\ref{app:proofs}.
For sufficiently large fixed
constants $A_0,A_1>0$, define
\[
e_n
=
2r^{m_n}
+
A_0\sqrt{
\frac{V_{m_n}\log n}{n}
},
\qquad
h_n
=
A_1\sqrt{\frac{\log n}{n}}.
\]

The adaptive estimator is defined
by the comparison rule
\begin{equation}\label{eq:adaptive-selector}
\widehat\theta_n^{\mathrm{ad}}
=
\begin{cases}
\widehat\theta_{R,n},
&
|\widetilde\theta_{P,n}
-\widehat\theta_{R,n}|
\leq e_n+h_n,\\
\widetilde\theta_{P,n},
&
\text{otherwise}.
\end{cases}
\end{equation}

The comparison is motivated by
the different statistical behavior
of the two estimators.
On the covariance-null submodel,
both estimators target $\theta(Q)$,
but the ratio estimator has
parametric risk.
Over the unrestricted model,
the polynomial estimator remains
uniformly consistent at the
logarithmic minimax scale,
whereas the ratio estimator may
have nonvanishing asymptotic bias.

The selection rule is therefore
designed to favor the ratio
estimator when the two estimators
are sufficiently close and to
retain the polynomial estimator
otherwise.

The theoretical analysis requires
more than consistency of the
selection rule. On the covariance-null
submodel, the probability of
selecting the polynomial estimator
must be sufficiently small to
preserve the $n^{-1}$ risk.
On the unrestricted class,
the comparison threshold must
prevent the selected estimator
from deviating substantially
from the robust polynomial estimator.

These properties are established
using concentration inequalities
for the ratio estimator and the
median-of-means polynomial estimator.
The tuning constants above are
theoretical quantities; their
finite-sample calibration is
considered separately.

\begin{theorem} 
\label{thm:adaptive}
For fixed $0<a<b<\infty$,
there exist deterministic sequences
$e_n,h_n>0$ and an estimator
$\widehat\theta_n^{\mathrm{ad}}$,
depending only on the observations
and the known bounds $(a,b)$,
such that
\begin{equation}\label{eq:adaptive-risk}
\begin{aligned}
\sup_{Q\in\mathcal Q_0(a,b)}
\mathbb E_Q
\left[
(\widehat\theta_n^{\mathrm{ad}}
-\theta(Q))^2
\right]
&=O(n^{-1}),\\
\liminf_{n\to\infty}
\frac{
-\log\displaystyle\sup_{Q\in\mathcal Q(a,b)}
\mathbb E_Q
\left[
(\widehat\theta_n^{\mathrm{ad}}
-\theta(Q))^2
\right]
}{L_n}
&\geq2\lambda.
\end{aligned}
\end{equation}
In particular, the estimator
achieves parametric squared-error
risk on $\mathcal Q_0(a,b)$ while
attaining the optimal full-model
logarithmic minimax exponent
from Theorem~\ref{thm:minimax-rate}.
\end{theorem}

\noindent\emph{Proof.}
See Appendix~\ref{app:proofs}.
\begin{remark} 
\label{rem:adaptive-scope}
The adaptation result concerns
squared-error risk on the
covariance-null submodel and
the logarithmic minimax exponent
over the unrestricted class.
It does not imply that the
adaptive estimator is
asymptotically efficient on
$\mathcal Q_0(a,b)$, nor that
it adapts to every intermediate
smoothness or dependence class.

Furthermore, the selection rule
should not be interpreted as a
uniformly consistent test of
membership in $\mathcal Q_0(a,b)$.
The median-of-means construction
is used to control estimation
risk through deviation bounds,
rather than to guarantee
correct model classification.

The constants defining the
theoretical thresholds may
be conservative at moderate
sample sizes. Consequently,
empirically calibrated
implementations must be
distinguished from the procedure
covered by the uniform
theoretical guarantees.
\end{remark}

\section{Numerical experiments}\label{sec:sim}

We examine the finite-sample implications of the three main theoretical results:
polynomial estimation of the latent mean composition (Theorem~\ref{thm:minimax-rate}),
detection of nonzero latent covariance (Theorem~\ref{thm:testing}), and
adaptation to the covariance-null submodel (Theorem~\ref{thm:adaptive}).
We also assess sensitivity to polynomial degree and support specification.
These experiments illustrate the mechanisms underlying the theory at selected
mixing distributions; they do not establish uniform minimax claims or verify
sharp logarithmic exponents.

\subsection{Simulation designs and implementation}
\label{subsec:sim-design}

For each replication, independent latent pairs $(S_i,U_i)$
are generated, followed by conditionally independent
Poisson counts as in~\eqref{eq:conditional-poisson}.
Unless specified otherwise, $(a,b)=(1,5)$ and $u_0=0.4$;
the estimators are supplied with the stated support bounds.
The independent beta design uses $(\alpha,\beta)=(2,3)$,
while the linear and nonlinear dependence designs use
$\gamma=0.6$ and $\delta=1.2$, respectively.
These parameter values are fixed throughout the main
simulation experiment. We consider four designs:

\begin{enumerate}

\item \emph{Constant composition:}
$S\sim\mathrm{Unif}[a,b]$ and $U\equiv u_0$;
thus $\theta=u_0=0.4$ and
$\operatorname{Cov}(S,U)=0$.

\item \emph{Independent heterogeneous composition:}
$S\sim\mathrm{Unif}[a,b]$ and
$U\sim\mathrm{Beta}(\alpha,\beta)$ independently,
with $(\alpha,\beta)=(2,3)$.
Thus $\theta=\alpha/(\alpha+\beta)=0.4$ and
$\operatorname{Cov}(S,U)=0$.

\item \emph{Intensity-dependent composition:}
$S\sim\mathrm{Unif}[a,b]$ and
\[
U=u_0+\gamma\frac{S-(a+b)/2}{b-a},
\]
where $|\gamma|\leq 2\min(u_0,1-u_0)$.
For $\gamma=0.6$, the resulting composition
satisfies $U\in[0.1,0.7]$.
Here $\theta=u_0=0.4$,
\[
\operatorname{Cov}(S,U)
=\frac{\gamma(b-a)}{12}=0.2,
\qquad
\theta_R-\theta
=\frac{\gamma(b-a)}{6(a+b)}
=\frac{1}{15}.
\]

\item \emph{Nonlinear dependence with zero covariance:}
with $V=(S-a)/(b-a)\sim\mathrm{Unif}[0,1]$, set
\[
U=u_0+\delta\left\{(V-1/2)^2-\frac{1}{12}\right\},
\]
with $\delta=1.2$.
For the chosen parameter values,
$U\in[0.3,0.6]$, so the composition remains
within its admissible range.
Moreover, $\theta=u_0=0.4$ and
$\operatorname{Cov}(S,U)=0$, despite the
nonlinear dependence between $S$ and $U$.

\end{enumerate}

The main experiment uses
$n\in\{250,500,1000,2000,4000,10000\}$
and $B=1000$ replications per configuration.
All methods are compared on the same generated dataset
within each replication.
We report the empirical ratio, projected polynomial
estimators of degrees $m=1,\ldots,8$, block-median
polynomial estimates, and theoretical and practically
calibrated adaptive rules.
The polynomial score is evaluated in the
falling-factorial basis,
\[
Z_{i,m}
=X_{i1}\sum_{j=0}^{m-1}
c_{j,m}(D_i-1)_j\,\mathbf1\{D_i\geq1\},
\]
not by substituting $D_i-1$ into the ordinary polynomial.

The ordinary polynomial estimator averages these scores
and is projected onto $[0,1]$ when required.
The robust version takes a median of block means.
The theoretical selector uses the degree sequence,
number of blocks and conservative comparison thresholds
specified in Section~\ref{subsec:adaptive}.
In contrast, the practical selector uses a fixed
polynomial degree $m=3$ and empirically calibrated
block and threshold choices.
We denote the resulting practically calibrated
adaptive estimator by
$\widehat{\theta}_n^{\mathrm{ad,pr}}$.
This estimator is distinct from the theoretical
adaptive estimator $\widehat{\theta}_n^{\mathrm{ad}}$
defined in Section~\ref{subsec:adaptive}
and is not covered by Theorem~\ref{thm:adaptive}.

The empirical MSE of an estimator $\widehat\theta$ is
\[
\frac{1}{B}\sum_{b=1}^{B}
\left(\widehat\theta^{(b)}-\theta\right)^2.
\]
We also report Monte Carlo bias, selection frequencies
and, for the comparisons in
Figure~\ref{fig:adaptive-selection},
95\% Monte Carlo uncertainty intervals.
Within each configuration, comparisons of estimators
use the same generated datasets.
No numerical instability was flagged in the recorded
main-campaign summaries; this does not constitute
a general numerical stability guarantee.

\subsection{Estimation accuracy and adaptive selection}\label{subsec:sim-accuracy}
\begin{figure}[!tbp]
\centering\includegraphics[width=\linewidth]{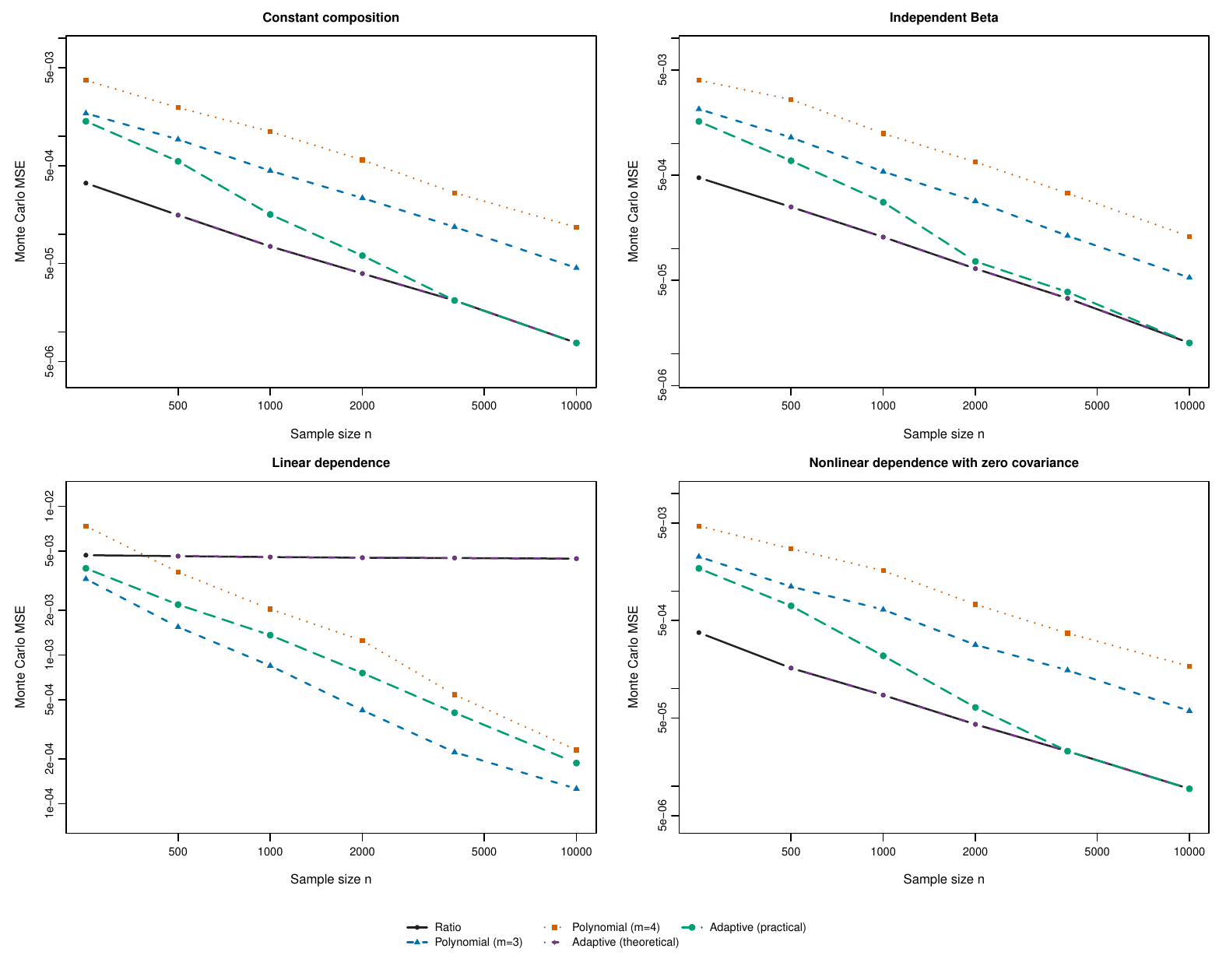}
\caption{Monte Carlo mean squared error as a function of sample size across the four latent-composition designs ($B=1000$ replications per configuration). The comparison includes the ratio estimator, polynomial estimators of degrees three and four, and the theoretical and practically calibrated adaptive estimators. Vertical axes use logarithmic scales. The ratio retains a non-vanishing error under linear dependence, whereas the polynomial estimators substantially reduce this error. The practical adaptive estimator is empirically tuned and is not covered by Theorem~\ref{thm:adaptive}.}\label{fig:mse-n}
\end{figure}

Table~\ref{tab:sim-main-mse} reports empirical mean squared errors at the largest sample size, (n=10,000), to provide a concise numerical comparison of the ratio, degree-three polynomial, and practically calibrated adaptive estimators. Results for all six sample sizes are displayed in Figures~\ref{fig:mse-n}–4, which document the finite-sample evolution of estimation error and adaptive selection.

The large-sample comparison highlights the persistent structural bias of the ratio estimator under nonzero latent covariance, in contrast to its favorable performance under covariance nullity. In the three covariance-null designs, including the nonlinear dependent design, the ratio estimator accurately targets $\theta(Q)$. Under linear dependence, its nonvanishing asymptotic bias dominates its MSE. The degree-three polynomial estimator substantially
reduces this error, whereas
$\widehat{\theta}_n^{\mathrm{ad,pr}}$
retains much of the ratio estimator's advantage
under zero covariance and performs substantially
better than the ratio under linear dependence. The latter improvement does not imply superiority over the ordinary polynomial estimator.

\begin{table}[!tbp]
\caption{Empirical MSE at $n=10\,000$ ($B=1000$). Values are rounded; ``practical'' denotes the empirically calibrated adaptive rule.}\label{tab:sim-main-mse}
\centering\small
\begin{tabular}{lrrr}
\hline
Design
& $\widehat\theta_{R,n}$
& $\widehat\theta_{n,3}$
& $\widehat\theta_n^{\mathrm{ad,pr}}$ \\
\hline
Constant composition & 0.000008 & 0.000045 & 0.000008\\
Independent beta & 0.000013 & 0.000053 & 0.000013\\
Linear dependence & 0.004450 & 0.000126 & 0.000187\\
Nonlinear, zero covariance & 0.000009 & 0.000059 & 0.000009\\
\hline
\end{tabular}
\end{table}

Figure~\ref{fig:mse-n} illustrates the distinction between estimating
$\theta(Q)$ and estimating the intensity-weighted ratio $\theta_R(Q)$.
The persistent ratio error under linear dependence agrees with
\eqref{eq:ratio-bias}; the low-degree polynomial estimators reduce this
structural bias at a variance cost. These fixed-degree results do not
validate the degree sequence or the sharp minimax exponent of
Theorem~\ref{thm:minimax-rate}.

\begin{figure}[!tbp]
\centering\includegraphics[width=\linewidth]{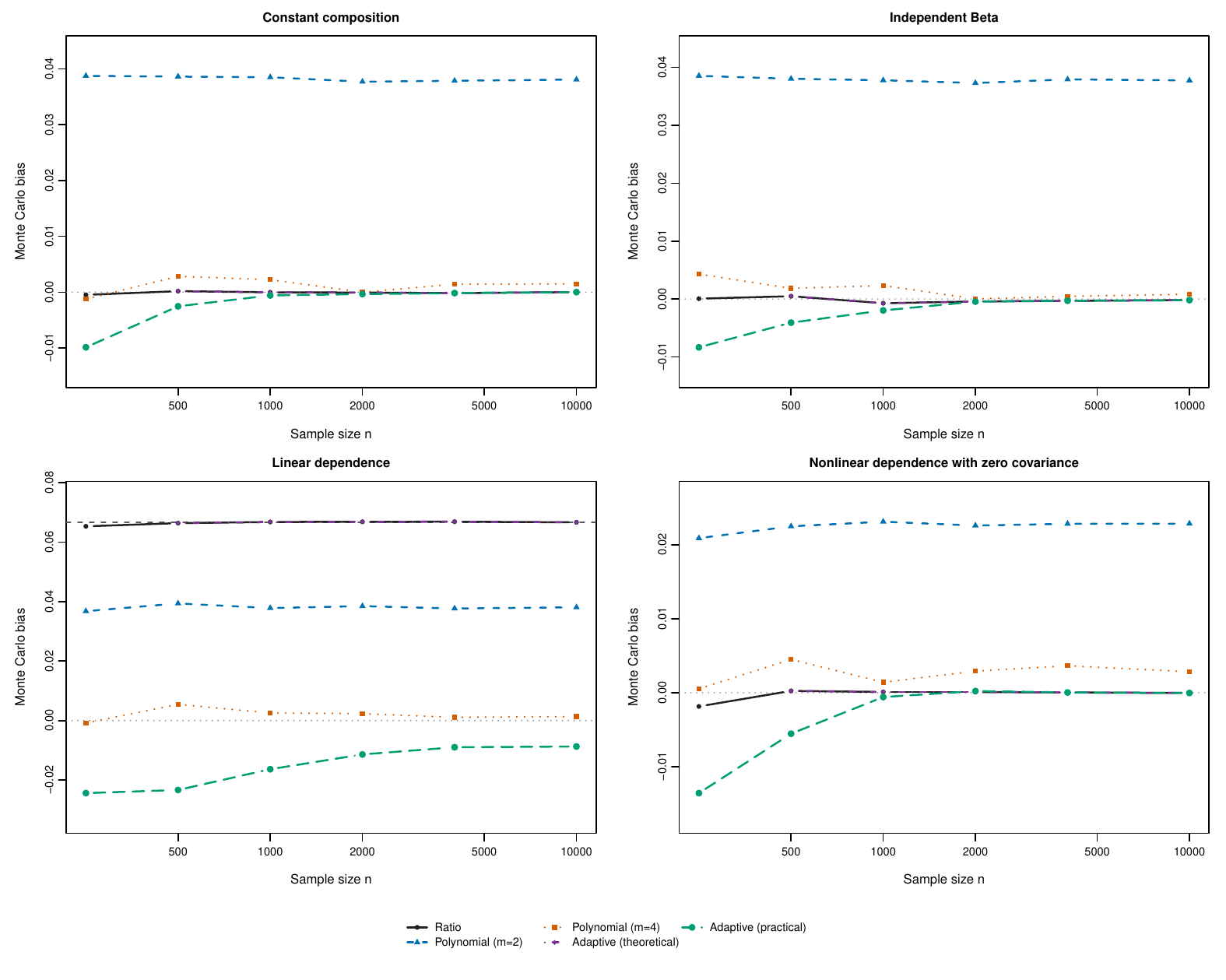}
\caption{Monte Carlo bias as a function of sample size in the four simulation designs ($B=1000$). The degree-two and degree-four polynomial estimators are displayed alongside the ratio and adaptive estimators. Under linear dependence, the ratio exhibits a persistent positive bias consistent with $\operatorname{Cov}(S,U)/\mathbb ES$; the low-degree polynomial estimators exhibit approximation bias. The polynomial degrees in this figure differ from those in Figure~\ref{fig:mse-n}.}\label{fig:bias-n}
\end{figure}

Figure~\ref{fig:bias-n} separates approximation bias from the structural
bias of the ratio. For $(a,b)=(1,5)$, the latter equals
$\theta_R-\theta=\gamma(b-a)/{6(a+b)}=\gamma/9$.
The remaining polynomial bias is compatible with the finite-degree
approximation error in Section~\ref{sec:minimax}.

\begin{figure}[!tbp]
\centering\includegraphics[width=\linewidth]{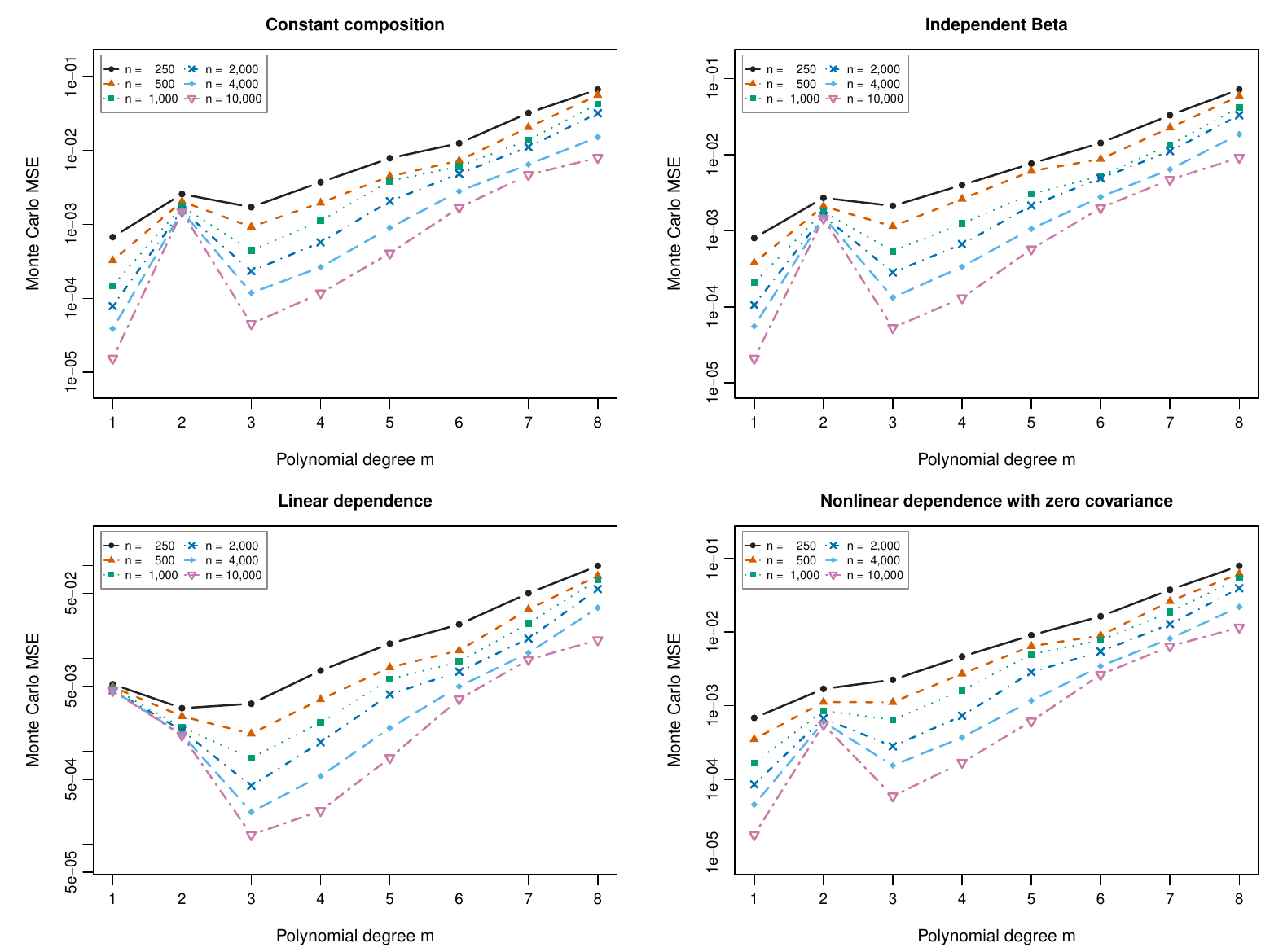}
\caption{Monte Carlo mean squared error of the polynomial estimator versus polynomial degree $m=1,\ldots,8$ at six sample sizes ($B=1000$). The logarithmic vertical scale highlights the bias--variance trade-off: increasing the degree initially improves accuracy in some designs, especially under linear dependence, but high degrees lead to markedly larger variance and MSE. The optimal finite-sample degree depends on the design and sample size.}\label{fig:degree}
\end{figure}

Figure~\ref{fig:degree} highlights the competing terms in the risk bound
of Section~\ref{sec:minimax}: geometric approximation error and increasing
variance of high-order factorial-moment scores. Hence minimizing the
Chebyshev approximation error alone is not a finite-sample selection rule.

At $n=10\,000$ under linear dependence, the paired Monte Carlo MSE
difference between the ratio and the ordinary degree-three polynomial
estimator is $0.004324$ (Monte Carlo standard error $0.000016$), with an
approximate 95\% Monte Carlo interval $[0.004293,0.004356]$.
The adaptive-minus-polynomial difference is approximately $0.000062$
(Monte Carlo standard error $0.000006$). This latter difference combines
the effects of median-of-means robustification and adaptive selection;
we separate them below.

\begin{figure}[!tbp]
\centering
\includegraphics[width=\linewidth]{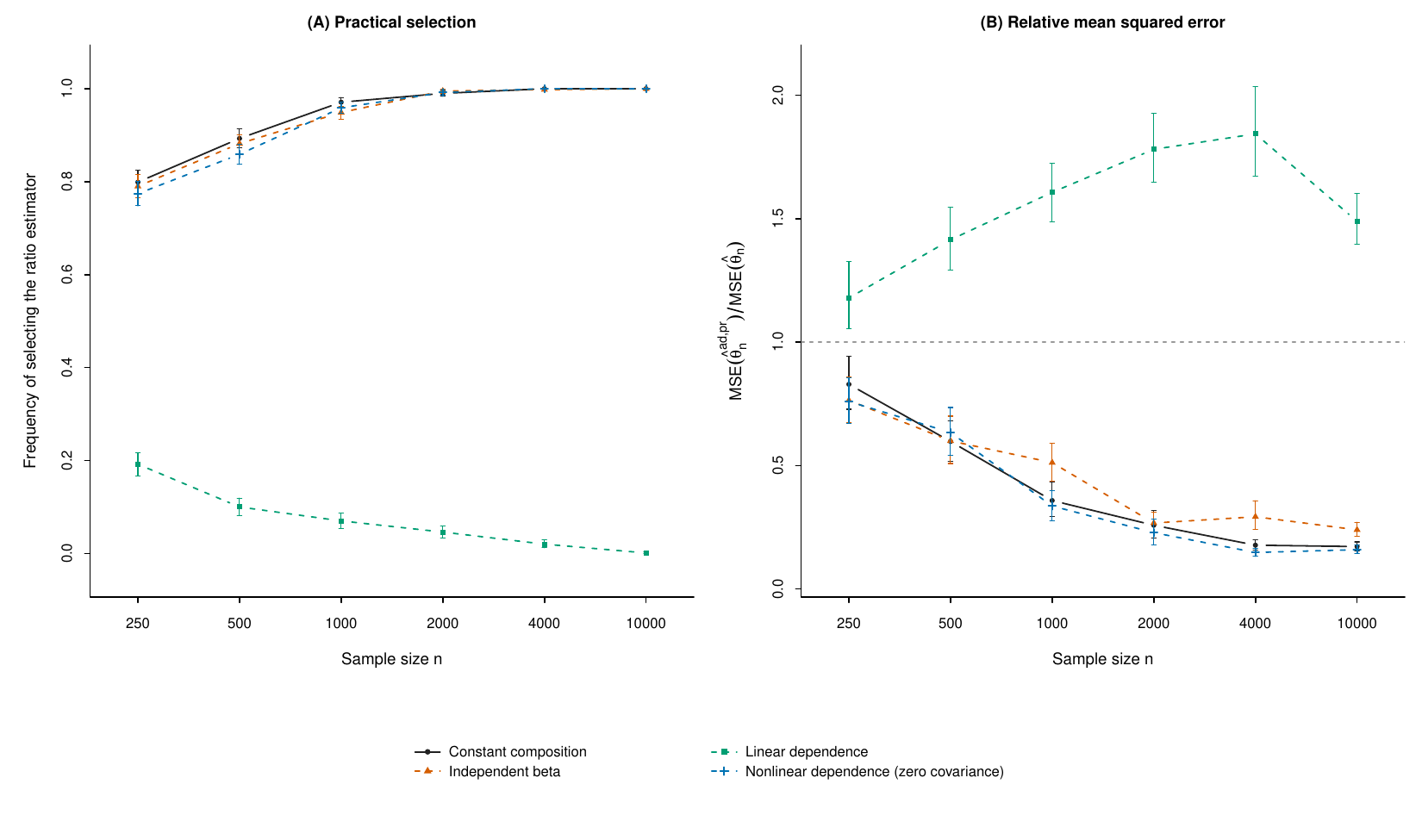}
\caption{Practical adaptive selection and relative mean squared error
($B=1000$ Monte Carlo replications per configuration).
Panel~(A) displays the frequency of selecting the ratio estimator.
Panel~(B) reports the ratio of the MSE of the practically calibrated
adaptive estimator to that of the ordinary degree-three polynomial
estimator; the horizontal dashed line corresponds to equal MSE.
Vertical bars are 95\% Monte Carlo percentile intervals based on
2,000 bootstrap resamples; in panel~(B), resampling is paired across
estimators within each replication. The practical rule is a finite-sample
implementation and does not inherit the uniform guarantees of
Theorem~\ref{thm:adaptive}.}
\label{fig:adaptive-selection}
\end{figure}

Figure~\ref{fig:adaptive-selection} shows the practical benefit of
favoring the ratio under the three covariance-null designs. At
$n=10\,000$, adaptive-to-polynomial MSE ratios are approximately
$0.17$, $0.24$, and $0.16$, respectively. The nonlinear design is
particularly informative: $S$ and $U$ are dependent although their
covariance vanishes. These results illustrate the structural submodel in
Theorem~\ref{thm:adaptive}, but do not establish its uniform risk bound.

Under linear dependence, ratio selection falls from $19.2\%$ at
$n=250$ to $0.1\%$ at $n=10\,000$, yet adaptive-to-polynomial MSE
ratios exceed one throughout (up to $1.84$ at $n=4\,000$).
Thus a low selection frequency need not imply a small risk penalty
when the ratio has a nonzero population bias.

To distinguish the selection cost from the cost of robustification, let
$A_n$ denote the event that the practical rule selects the ratio and let
$\widetilde\theta_{P,n}$ denote the median-of-means polynomial estimator
used by that rule. For the same data and the same median-of-means
realization, the adaptive estimate equals
$\widehat\theta_{R,n}\mathbf1_{A_n}+
\widetilde\theta_{P,n}\mathbf1_{A_n^c}$. Consequently,
\[
\begin{aligned}
&\mathbb E_Q[(\widehat\theta_{\mathrm{ad},n}-\theta)^2]
 -\mathbb E_Q[(\widetilde\theta_{P,n}-\theta)^2]\\
&\qquad=\mathbb E_Q\!\left[
\big\{(\widehat\theta_{R,n}-\theta)^2
 -(\widetilde\theta_{P,n}-\theta)^2\big\}
\mathbf1_{A_n}\right].
\end{aligned}
\]
The empirical analogue is computed by pairing the three methods within
each Monte Carlo replication. Table~\ref{tab:sim-selection-cost} reports
this diagnostic for the nonzero-covariance design, where the trade-off is
most pronounced.

\begin{table}[!tbp]
\caption{Finite-sample cost of practical adaptive selection
under linear dependence ($B=1000$). The relative selection
cost is the paired MSE difference between
$\widehat\theta_n^{\mathrm{ad,pr}}$ and the
median-of-means polynomial estimator, divided by the
MSE of the latter. The last column reports the MSE
ratio of $\widehat\theta_n^{\mathrm{ad,pr}}$ to the
ordinary degree-three polynomial estimator
$\widehat\theta_{n,3}$. Values are rounded.}
\label{tab:sim-selection-cost}
\centering\small
\begin{tabular}{rrrr}
\hline
$n$ &
Ratio selection (\%) &
Relative selection cost (\%) &
$\dfrac{\mathrm{MSE}(\widehat\theta_n^{\mathrm{ad,pr}})}
{\mathrm{MSE}(\widehat\theta_{n,3})}$ \\
\hline
250   & 19.2 & 19.7 & 1.18\\
500   & 10.0 & 14.8 & 1.41\\
1000  &  7.0 & 21.4 & 1.61\\
2000  &  4.6 & 25.6 & 1.78\\
4000  &  2.0 & 23.1 & 1.84\\
10000 &  0.1 &  2.1 & 1.49\\
\hline
\end{tabular}
\end{table}

Table~\ref{tab:sim-selection-cost} separates selection from
robustification. At $n=2\,000$, selection increases MSE by $25.6\%$
relative to the median-of-means polynomial estimator, despite selecting
the ratio in only $4.6\%$ of replications. At $n=10\,000$, the
corresponding selection cost is $2.1\%$; the ordinary polynomial,
median-of-means polynomial, and practical adaptive MSEs are approximately
$1.26\times10^{-4}$, $1.84\times10^{-4}$, and $1.87\times10^{-4}$.
Thus robustification accounts for most of the observed excess over the
ordinary polynomial at that sample size. The $0.1\%$ selection frequency
represents a single event in 1,000 replications, so its risk contribution
is estimated imprecisely.

The theoretical selector of Section~\ref{subsec:adaptive} has a different finite-sample
behavior. Its conservative threshold leads to ratio selection in all
reported replications at $n=10\,000$ under linear dependence, yielding
MSE $0.004450$, compared with $0.000187$ for the practical selector.
This does not contradict Theorem~\ref{thm:adaptive}, which concerns
asymptotic uniform risk bounds and does not guarantee effective selection
at these sample sizes. Conversely, the practical fixed-degree rule is not
a numerical verification of that theorem.

\subsection{Covariance testing, local alternatives, and sensitivity analyses}\label{subsec:sim-testing}

The covariance test is calibrated by an independent simulation under a specified nonlinear covariance-zero design, using $B_{\mathrm{cal}}=1500$ calibration replicates. We also evaluate empirical rejection under constant-composition and independent-beta null designs. This procedure is a finite-sample diagnostic, not a proof of uniform level over $\mathcal Q_0(a,b)$.

\begin{figure}[!tbp]
\centering\includegraphics[width=\linewidth]{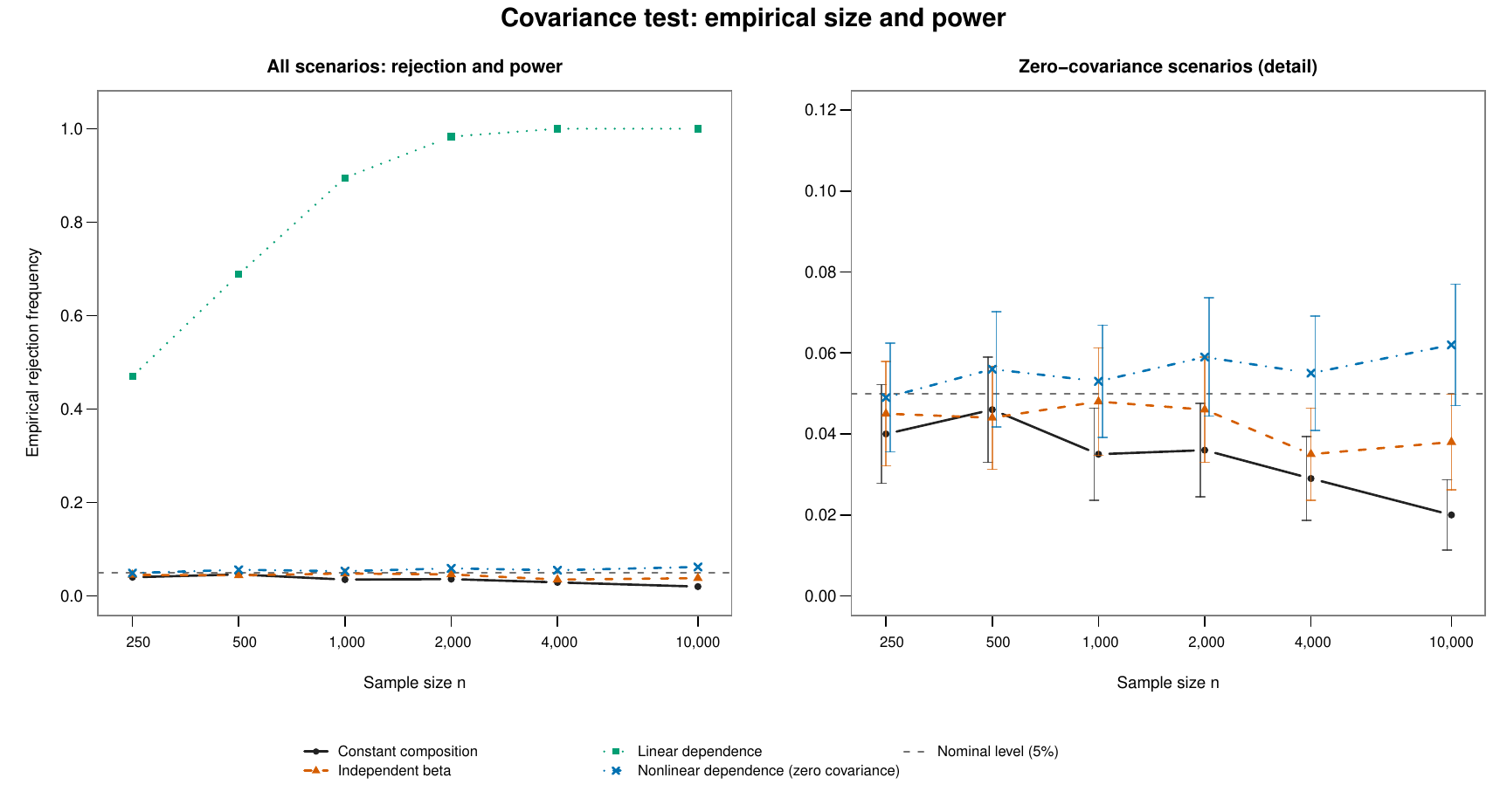}
\caption{Covariance test: empirical rejection frequencies and power ($B=1000$ evaluation replications per configuration). Left: all four simulation designs, including the fixed linear-dependence alternative. Right: magnified view of the three zero-covariance designs, with approximate 95\% Monte Carlo confidence intervals $\widehat p\pm1.96\sqrt{\widehat p(1-\widehat p)/B}$. The dashed horizontal line denotes the nominal 5\% level. Critical values were calibrated using a separate sample of $B_{\mathrm{cal}}=1500$ replications under a specified nonlinear zero-covariance design; these plots do not establish uniform size control over $\mathcal Q_0(a,b)$.}\label{fig:cov-test}
\end{figure}

Figure~\ref{fig:cov-test} shows increasing rejection probability against the fixed linear-dependence alternative, approaching one at larger sample sizes. Under the three covariance-null designs, empirical rejection frequencies remain of the same order as the nominal level, but are not identical: the constant-composition case is conservative at larger $n$, whereas the nonlinear zero-covariance case is slightly above 5\% in some configurations. The displayed error bars quantify Monte Carlo uncertainty conditional on the implemented calibration procedure, not uncertainty due to the independent critical-value calibration.

To examine local alternatives, we use the linear design with $\gamma_n$ decreasing in $n$. Write
\[
\lambda=\log\!\left(\frac{\sqrt b+\sqrt a}{\sqrt b-\sqrt a}\right),\qquad
\rho_n(\kappa)=\exp\!\left\{-\kappa\lambda\frac{\log n}{\log\log n}\right\},\qquad
\gamma_n=\frac{12\rho_n(\kappa)}{b-a}.
\]
The admissibility condition on $\gamma_n$ is checked before generating data. 

Table~\ref{tab:sim-local} reports the empirical
rejection frequencies under local linear alternatives.
The polynomial degree is fixed at $m=3$ for
$n\leq 4000$ and increases to $m=4$ at $n=10\,000$.
Consequently, the last row reflects both an increase
in sample size and a change in polynomial degree,
and cannot be interpreted as a pure sample-size effect.
A complementary fixed-sample-size experiment,
reported below, examines the sensitivity of
empirical power to this degree change.

\begin{table}[!tbp]
\caption{Empirical rejection frequencies under local
linear alternatives ($B=1000$). The latent covariance is
$\rho_n=\exp\{-\kappa\lambda\log(n)/\log\log(n)\}$,
where
$\lambda=\log\{(\sqrt{b}+\sqrt{a})/
(\sqrt{b}-\sqrt{a})\}$.
Critical values are calibrated using
$B_{\mathrm{cal}}=1500$ independent null replications
at nominal level $\alpha=0.05$.
The polynomial degree is
$m=\lfloor\log(n)/\log\log(n)\rfloor$.}
\label{tab:sim-local}
\centering\small
\begin{tabular}{rrrrrr}
\hline
$n$ & $m$ & $\kappa=0.5$ & $\kappa=0.75$
& $\kappa=1$ & $\kappa=1.25$\\
\hline
250   & 3 & 0.488 & 0.105 & 0.061 & 0.033\\
500   & 3 & 0.707 & 0.187 & 0.070 & 0.062\\
1000  & 3 & 0.810 & 0.251 & 0.065 & 0.039\\
2000  & 3 & 0.966 & 0.359 & 0.075 & 0.036\\
4000  & 3 & 0.992 & 0.517 & 0.103 & 0.046\\
10000 & 4 & 0.953 & 0.191 & 0.036 & 0.022\\
\hline
\end{tabular}
\end{table}

The ordering of rejection frequencies across $\kappa$
is consistent with a challenging detection problem,
although these linear alternatives need not be least
favorable for the unrestricted testing problem.
The decrease in rejection frequencies at $n=10\,000$
must be interpreted cautiously, since both the local
alternative and the polynomial degree change with
sample size.

To assess the role of the polynomial degree, we
conducted an additional experiment at $n=10\,000$,
comparing $m=3$ and $m=4$ under identical local
alternatives, using paired Monte Carlo samples and
separately calibrated critical values.
For $\kappa=0.75$, the empirical rejection frequencies
were $0.755$ and $0.222$, respectively, yielding a
paired difference of $0.533$ (95\% Monte Carlo
confidence interval: $[0.491,0.575]$).
Independent null simulations yielded rejection
frequencies of $0.0455$ and $0.0525$, respectively,
close to the nominal level of $0.05$.

These findings demonstrate substantial sensitivity
of finite-sample power to the polynomial degree
in this configuration. They neither establish
uniform size control over the covariance-null class
nor validate the sharp logarithmic separation
exponent in Theorem~\ref{thm:testing}.

\paragraph{Sensitivity to polynomial degree and support bounds.}\label{subsec:sim-support}

A controlled experiment at $n=10\,000$ compares degrees three and four under the same local-alternative family, with separate critical-value calibration for each degree. The calibrated critical values are $0.041992$ and $0.074407$, respectively. At $\kappa=0.75$, empirical power falls from $0.752$ at degree three to $0.207$ at degree four; at $\kappa=1$, it falls from $0.128$ to $0.046$. This demonstrates a substantial finite-sample bias--variance trade-off in the polynomial regularization. Under the three investigated null laws, rejection frequencies are at most $0.050$ for both degrees, but this remains a pointwise calibration check.

To assess support sensitivity, we compare the correctly specified intervals $[2,5]$, $[1,5]$, and $[0.2,5]$, and examine deliberately misspecified working intervals while keeping the reference data-generating support fixed. Table~\ref{tab:sim-support} reports degree-three polynomial MSE at $n=10\,000$. The near-zero lower bound markedly increases MSE, especially under linear dependence. A misspecified narrower working interval can produce a smaller MSE in this particular configuration, but it violates the support assumptions used to justify the polynomial approximation and should not be interpreted as a generally preferable choice.

\begin{table}[!tbp]
\caption{Degree-three polynomial MSE at $n=10\,000$ under different true or working support bounds ($B=1000$). The last two rows change working bounds only.}\label{tab:sim-support}
\centering\small
\begin{tabular}{lrr}
\hline
Support configuration & Constant composition & Linear dependence\\
\hline
Correct $[2,5]$ & 0.000017 & 0.000027\\
Correct $[1,5]$ & 0.000049 & 0.000122\\
Correct $[0.2,5]$ & 0.000330 & 0.001723\\
Working $[1.5,4.5]$ & 0.000032 & 0.000105\\
Working $[0.2,6]$ & 0.000223 & 0.004669\\
\hline
\end{tabular}
\end{table}

Finally, the sign experiment uses $\gamma=\pm0.6$. Since $(a,b)=(1,5)$, the ratio's theoretical asymptotic bias is $\gamma/9=\pm0.066667$. At $n=10\,000$, the empirical biases are $-0.066626$ and $0.066715$, respectively. These values were calculated after grouping simulation replications by the sign of $\gamma$; pooling the two signs would spuriously cancel their biases.

\paragraph{Interpretation and limitations.}\label{subsec:sim-discussion}

The numerical findings support three qualitative conclusions.
First, covariance nullity, not independence, determines whether the
ratio targets $\theta(Q)$. Second, polynomial approximation trades
structural bias for degree-dependent variance, with substantial
sensitivity to the lower support bound. Third, practical adaptation
may reduce risk on selected covariance-null laws but incur costs from
both robustification and occasional ratio selection under alternatives.
These are pointwise, finite-sample observations: the sharp minimax
exponent of Theorem~\ref{thm:minimax-rate}, the uniform separation
boundary of Theorem~\ref{thm:testing}, and the uniform adaptation
properties of Theorem~\ref{thm:adaptive} rest on their proofs, not on
simulation evidence.

\section{Discussion}\label{sec:discussion}

This paper has investigated the estimation of mean latent
composition in bivariate Poisson mixtures with an unknown
joint mixing distribution. We established a sharp logarithmic
minimax exponent, revealing that identifiability does not
necessarily imply accurate estimation. The matching upper
and lower bounds show that the slow convergence rate is
intrinsic to the statistical model, reflecting the difficulty
of recovering latent information from observed counts.
We also established a sharp logarithmic detection boundary
for testing zero latent covariance and showed that this
restriction, which is weaker than independence, permits
estimation at the parametric rate. An adaptive procedure
achieves this faster rate on the covariance-null submodel
while retaining the optimal logarithmic minimax exponent
over the unrestricted class. The numerical experiments
illustrate these theoretical findings and highlight the
finite-sample costs associated with polynomial regularization
and adaptive selection.
The analysis assumes known and fixed bounds on the latent
total intensity, with a strictly positive lower bound.
Extending the results to unknown support bounds, mixing
distributions allowing intensities near zero, and additional
classes of latent dependence remains an important direction
for future research. Developing fully data-driven
regularization and adaptive procedures with theoretical
guarantees is another natural extension.

\appendix
\section{Proofs of the main results and auxiliary bounds}\label{app:proofs}
Throughout this appendix, $0<a<b<\infty$ are fixed, constants may depend on $(a,b)$ only, $D=X_1+X_2$, and $P_Q$ is the law of one observed pair. We write $\|\nu\|_{\mathrm{var}}=|\nu|([a,b])$ for the total variation \emph{norm} of a signed measure, and $d_{\mathrm{TV}}(P,P')=\sup_A|P(A)-P'(A)|=\frac12\sum_y|P(y)-P'(y)|$ for the total variation \emph{distance} between probability laws.

\subsection*{Proof of Theorem~\ref{thm:identifiability}}
\begin{proof}
Let $G=Q\circ\Phi^{-1}$ and $G'=Q'\circ\Phi^{-1}$ be the laws of $(SU,S(1-U))$ on the compact set $\mathcal L_{a,b}$. For all integers $j,k\ge0$, conditional independence and Poisson factorial moments give
\[
\mathbb E_Q[(X_1)_j(X_2)_k]=\int\lambda_1^j\lambda_2^k\,G(d\lambda_1,d\lambda_2).
\]
Equality of observable laws implies equality of integrals of all bivariate polynomials. These polynomials form a unital algebra separating points of $\mathcal L_{a,b}$; Stone--Weierstrass yields density in $C(\mathcal L_{a,b})$. Thus $G=G'$, and the continuous inverse $\Phi^{-1}(\lambda_1,\lambda_2)=(\lambda_1+\lambda_2,\lambda_1/(\lambda_1+\lambda_2))$ gives $Q=Q'$.
\end{proof}

\begin{lemma}[Uniform upper bound]\label{lem:minimax-upper}
There exist $C_0,C_1\in(0,\infty)$ depending only on $(a,b)$ such that, for every $m,n\ge1$,
\begin{equation}\label{eq:upper-bound}
\sup_{Q\in\mathcal Q(a,b)}\mathbb E_Q[(\widehat\theta_{n,m}-\theta(Q))^2]\le4r^{2m}+\frac{C_0C_1^m m!}{n}.
\end{equation}
Moreover $\sup_Q\mathbb E_Q[Z_{i,m}^2]\le C_0C_1^m m!$.
\end{lemma}
\begin{proof}
The Chebyshev approximation bound in Section~\ref{subsec:minimax-upper} implies $|\mathbb E_Q Z_{i,m}-\theta(Q)|\le2r^m$.
For $k\ge1$, the factorial-product identity and the Poisson factorial moments imply, conditional on $S=s$,
\[
\mathbb E[(D)_k^2\mid S=s]=\sum_{\ell=0}^k\binom{k}{\ell}^2\ell!s^{2k-\ell}
=k!s^k\sum_{h=0}^k\binom{k}{h}\frac{s^h}{h!}
\le e^b k!(2b)^k.
\]
Since $0\le H_j(Y)\le(D)_{j+1}$, this bounds $\|H_j\|_{L^2(Q)}$ uniformly.

For a polynomial $p$, write $\|p\|_{\mathrm{coef},1}$ for the sum of absolute values of its monomial coefficients. The recurrence $T_{m+1}=2xT_m-T_{m-1}$ yields $\|T_{m+1}\|_{\mathrm{coef},1}\le2\|T_m\|_{\mathrm{coef},1}+\|T_{m-1}\|_{\mathrm{coef},1}$. Since $\alpha=1+\sqrt2$ satisfies $\alpha^2=2\alpha+1$, induction from $T_0=1,T_1=x$ gives $\|T_m\|_{\mathrm{coef},1}\le\alpha^m\le3^m$.
Set $z=(a+b)/(b-a)$ and $\kappa=2/(b-a)$. Because $z+\kappa>1$, substitution gives
\[
\|T_m(z-\kappa s)\|_{\mathrm{coef},1}\le3^m(z+\kappa)^m.
\]
Since $T_m(z)\ge1$ and the constant coefficient of $P_m(s)=T_m(z-\kappa s)/T_m(z)$ is one, division of $1-P_m(s)$ by $s$ only shifts nonconstant coefficients. Consequently,
\[
\sum_{j=0}^{m-1}|c_{j,m}|\le A^m,\qquad A=3(z+\kappa).
\]
Put $B=\max\{1,2b\}$. For $j+1\le m$, $\|H_j\|_2\le\sqrt{e^b m!B^m}$. Minkowski's inequality now gives
\[
\|Z_{i,m}\|_2\le A^m\sqrt{e^b m!B^m}.
\]
Thus $\sup_Q\mathbb E_Q[Z_{i,m}^2]\le C_0C_1^m m!$ with the explicit choice
\[
C_0=e^b,\qquad C_1=9\left(\frac{a+b+2}{b-a}\right)^2\max\{1,2b\}.
\]
Independence and the squared bias bound finish the proof.
\end{proof}

\begin{lemma}[Two-point lower bound]
\label{lem:minimax-lower}
For every pair of integers $m,n\ge1$,
\begin{equation}\label{eq:lower-bound}
R_n^*(a,b)
\ge
\frac{1}{8T_m(z)^2}
\left[
1-\frac{n(2b)^{m+1}}{(m+1)!}
\right]_+,
\end{equation}
where $z=(a+b)/(b-a)>1$ and $[x]_+=\max\{x,0\}$.
\end{lemma}

\begin{proof}
Let
\[
V_m=\operatorname{span}\{s,s^2,\ldots,s^m\}
\subset C([a,b]).
\]
Since every polynomial in $1-V_m$ has degree at
most $m$ and takes the value $1$ at zero, the
Chebyshev extremal property for uniform
approximation (see \cite{DeVoreLorentz1993})
gives
\[
E_m:=
\operatorname{dist}_{\infty}(1,V_m)
=\frac{1}{T_m(z)}.
\]
Indeed, the polynomial
\[
P_m(s)=
\frac{
T_m\left(\frac{a+b-2s}{b-a}\right)
}{
T_m(z)
}
\]
satisfies $P_m(0)=1$ and attains alternating
extrema of magnitude $E_m$ at $m+1$ points
of $[a,b]$. The alternation theorem therefore
establishes its optimality.

Since $V_m$ is a closed subspace of
$C([a,b])$, the Hahn--Banach theorem and
the Riesz representation theorem yield a
finite signed Borel measure $\nu_m$ on
$[a,b]$ such that
\[
\|\nu_m\|_{\mathrm{var}}=1,
\qquad
\int s^j\,\nu_m(ds)=0,
\quad j=1,\ldots,m,
\]
and
\[
\int\nu_m(ds)=E_m,
\]
where the sign of $\nu_m$ has been chosen
so that the last integral is positive.
Here $\|\nu_m\|_{\mathrm{var}}$ denotes
the total variation norm of a signed measure.

Write the Jordan decomposition as
$\nu_m=\nu_m^+-\nu_m^-$.
Since
\[
\nu_m^+([a,b])+\nu_m^-([a,b])=1
\]
and
\[
\nu_m^+([a,b])-\nu_m^-([a,b])=E_m,
\]
we obtain
\[
\nu_m^+([a,b])=\frac{1+E_m}{2},
\qquad
\nu_m^-([a,b])=\frac{1-E_m}{2}.
\]
Consequently,
$F_m=\nu_m^++\nu_m^-$ is a probability
measure on $[a,b]$.

Define two probability distributions on
$[a,b]\times[0,1]$ by
\[
\begin{aligned}
Q_+(ds,du)
&=\nu_m^+(ds)\delta_1(du)
  +\nu_m^-(ds)\delta_0(du),\\
Q_-(ds,du)
&=\nu_m^+(ds)\delta_0(du)
  +\nu_m^-(ds)\delta_1(du).
\end{aligned}
\]
Both distributions belong to
$\mathcal Q(a,b)$ and have the same
$S$-marginal $F_m$. Moreover,
\[
\theta(Q_+)=\frac{1+E_m}{2},
\qquad
\theta(Q_-)=\frac{1-E_m}{2},
\]
so that
\begin{equation}\label{eq:two-point-gap}
\theta(Q_+)-\theta(Q_-)=E_m.
\end{equation}

Let $P_+$ and $P_-$ denote the distributions
of one observed count vector under
$Q_+$ and $Q_-$, respectively.
Since $U\in\{0,1\}$ under both mixing
distributions, the observed count vectors
are supported on the coordinate axes.

For every integer $d\ge1$, define
\[
a_d=\frac{1}{d!}
\int_{a}^{b}e^{-s}s^d\,\nu_m(ds).
\]
The signed probability mass differences
satisfy
\[
P_+(d,0)-P_-(d,0)=a_d,
\]
and
\[
P_+(0,d)-P_-(0,d)=-a_d.
\]
At $(0,0)$, the two probability masses
coincide because $Q_+$ and $Q_-$ have
the same $S$-marginal.

Using the convention
\[
d_{\mathrm{TV}}(P_+,P_-)
=\frac12\sum_{y\in\mathbb N_0^2}
|P_+(y)-P_-(y)|,
\]
we therefore obtain
\begin{equation}\label{eq:tv-axis}
d_{\mathrm{TV}}(P_+,P_-)
=\sum_{d=1}^{\infty}|a_d|.
\end{equation}

We first consider $1\le d\le m$.
Taylor's formula gives
\[
e^{-s}
=
\sum_{k=0}^{m-d}\frac{(-s)^k}{k!}
+R_{m-d+1}(s),
\]
where, for $s\ge0$,
\[
|R_{m-d+1}(s)|
\le
\frac{s^{m-d+1}}{(m-d+1)!}.
\]
By the moment cancellation conditions,
\[
\int s^{d+k}\,\nu_m(ds)=0,
\qquad 0\le k\le m-d.
\]
Consequently,
\[
a_d=
\frac1{d!}
\int s^dR_{m-d+1}(s)\,\nu_m(ds).
\]
Since $\|\nu_m\|_{\mathrm{var}}=1$
and $s\le b$ on its support,
\begin{equation}\label{eq:ad-bound}
|a_d|
\le
\frac{b^{m+1}}
{d!(m-d+1)!},
\qquad 1\le d\le m.
\end{equation}

For the remaining terms, the triangle
inequality and Tonelli's theorem yield
\[
\begin{aligned}
\sum_{d=m+1}^{\infty}|a_d|
&\le
\int
\left(
\sum_{d=m+1}^{\infty}
e^{-s}\frac{s^d}{d!}
\right)|\nu_m|(ds)\\
&=
\int
\Pr\{\operatorname{Pois}(s)\ge m+1\}
\,|\nu_m|(ds).
\end{aligned}
\]
If $N\sim\operatorname{Pois}(s)$, then
\[
\mathbf1_{\{N\ge m+1\}}
\le
\frac{(N)_{m+1}}{(m+1)!}.
\]
Taking expectations and using the
Poisson factorial-moment identity gives
\[
\Pr\{N\ge m+1\}
\le
\frac{s^{m+1}}{(m+1)!}.
\]
It follows that
\begin{equation}\label{eq:tail-bound}
\sum_{d=m+1}^{\infty}|a_d|
\le
\frac{b^{m+1}}{(m+1)!}.
\end{equation}

Combining
\eqref{eq:tv-axis},
\eqref{eq:ad-bound}, and
\eqref{eq:tail-bound}, we obtain
\[
\begin{aligned}
d_{\mathrm{TV}}(P_+,P_-)
&\le
\frac{b^{m+1}}{(m+1)!}
\left[
\sum_{d=1}^{m}
\binom{m+1}{d}+1
\right]\\
&=
\frac{(2^{m+1}-1)b^{m+1}}
{(m+1)!}\\
&\le
\frac{(2b)^{m+1}}{(m+1)!}.
\end{aligned}
\]
By the product coupling inequality,
\begin{equation}\label{eq:tv-product}
d_{\mathrm{TV}}
(P_+^{\otimes n},P_-^{\otimes n})
\le
n\,d_{\mathrm{TV}}(P_+,P_-).
\end{equation}

It remains to relate the testing
difficulty to the estimation risk.
Let $\widehat\theta$ be any estimator
based on $n$ independent observations,
and define the test
\[
\phi_n=
\mathbf1\left\{
\widehat\theta>
\frac{\theta(Q_+)+\theta(Q_-)}{2}
\right\}.
\]
By \eqref{eq:two-point-gap}, on
$\{\phi_n=1\}$ under $Q_-$,
\[
|\widehat\theta-\theta(Q_-)|
\ge\frac{E_m}{2},
\]
whereas on $\{\phi_n=0\}$ under $Q_+$,
\[
|\widehat\theta-\theta(Q_+)|
\ge\frac{E_m}{2}.
\]
Hence
\[
\begin{aligned}
&\mathbb E_{Q_-}
(\widehat\theta-\theta(Q_-))^2
+\mathbb E_{Q_+}
(\widehat\theta-\theta(Q_+))^2\\
&\quad\ge
\frac{E_m^2}{4}
\left[
P_-^{\otimes n}(\phi_n=1)
+
P_+^{\otimes n}(\phi_n=0)
\right].
\end{aligned}
\]
The sum of the two testing errors is
bounded below by
\[
1-d_{\mathrm{TV}}
(P_+^{\otimes n},P_-^{\otimes n}).
\]
Since the maximum of two nonnegative
numbers is at least their arithmetic
mean, the classical two-point argument
(see \cite[Chapter~2]{Tsybakov2009})
gives
\[
\begin{aligned}
R_n^*(a,b)
&\ge
\frac{E_m^2}{8}
\left[
1-d_{\mathrm{TV}}
(P_+^{\otimes n},P_-^{\otimes n})
\right]\\
&\ge
\frac{E_m^2}{8}
\left[
1-\frac{n(2b)^{m+1}}{(m+1)!}
\right]_+.
\end{aligned}
\]
Substituting $E_m=1/T_m(z)$
establishes \eqref{eq:lower-bound}.
\end{proof}

\begin{proof}[Proof of Theorem~\ref{thm:minimax-rate}]
Throughout the proof, let
\[
L_n=\frac{\log n}{\log\log n},
\qquad
r=\frac{\sqrt b-\sqrt a}{\sqrt b+\sqrt a},
\qquad
\lambda=-\log r>0.
\]
All asymptotic statements are understood as
$n\to\infty$, with $0<a<b<\infty$ fixed.

We first derive an upper bound on the minimax risk,
which yields the lower bound on its negative
logarithmic exponent. Set
$m_n=\lfloor L_n\rfloor$ for all sufficiently
large $n$. By Lemma~\ref{lem:minimax-upper},
\[
R_n^*(a,b)
\le
4r^{2m_n}
+
\frac{C_0C_1^{m_n}m_n!}{n},
\]
where $C_0,C_1>0$ depend only on $a$ and $b$.

Since $m_n=L_n+O(1)$, Stirling's formula gives
\[
\log(m_n!)
=
m_n\log m_n-m_n+O(\log m_n).
\]
Moreover, $\log m_n
=
\log\log n-\log\log\log n+o(1).$

Consequently,
\[
\begin{aligned}
\log\left(
\frac{C_1^{m_n}m_n!}{n}
\right)
&=
m_n\log C_1+m_n\log m_n
-m_n-\log n+O(\log m_n)\\
&=
-L_n\log\log\log n+O(L_n).
\end{aligned}
\]
As $\log\log\log n\to\infty$, it follows that
\[
\frac{1}{L_n}
\log\left(
\frac{C_1^{m_n}m_n!}{n}
\right)
\longrightarrow-\infty.
\]
On the other hand,
\[
\frac{\log(r^{2m_n})}{L_n}
=
-2\lambda\frac{m_n}{L_n}
\longrightarrow-2\lambda.
\]
Therefore,
\[
\frac{C_0C_1^{m_n}m_n!/n}
{r^{2m_n}}
\longrightarrow0,
\]
and hence $R_n^*(a,b)
\le
4r^{2m_n}\{1+o(1)\}.$
Taking negative logarithms and dividing
by $L_n$ yields
\begin{equation}\label{eq:proof-liminf}
\liminf_{n\to\infty}
\frac{-\log R_n^*(a,b)}{L_n}
\ge 2\lambda.
\end{equation}

We next derive a lower bound on the minimax risk,
which yields the reverse inequality for the
logarithmic exponent. Fix $\varepsilon>0$ and set $m_n=\lceil(1+\varepsilon)L_n\rceil.$
By Lemma~\ref{lem:minimax-lower},
\[
R_n^*(a,b)
\ge
\frac{1}{8T_{m_n}(z)^2}
[1-A_n]_+,
\]
where
\[
A_n=
\frac{n(2b)^{m_n+1}}{(m_n+1)!},
\qquad
z=\frac{a+b}{b-a}>1.
\]

Applying Stirling's formula,
\[
\begin{aligned}
\log A_n
&=
\log n+(m_n+1)\log(2b)\\
&\quad
-(m_n+1)\log(m_n+1)
+(m_n+1)
+O(\log(m_n+1)).
\end{aligned}
\]
Since $m_n=(1+\varepsilon)L_n+O(1)$
and $\log(m_n+1)
=
\log\log n-\log\log\log n
+\log(1+\varepsilon)+o(1),$
we obtain
\[
\log A_n
=
-\varepsilon\log n
+O(L_n\log\log\log n).
\]
Because $L_n\log\log\log n=o(\log n),$
it follows that $\log A_n
=
-\varepsilon\log n+o(\log n)
\longrightarrow-\infty.$
Thus $A_n\to0$, and consequently
\[
[1-A_n]_+=1-o(1).
\]

Since $0<r<1$, the Chebyshev identity gives
\[
\frac{1}{T_m(z)}
=
\frac{2r^m}{1+r^{2m}}
\ge r^m.
\]
Combining this inequality with
Lemma~\ref{lem:minimax-lower} and $A_n\to0$
yields
\[
R_n^*(a,b)
\ge
\frac18 r^{2m_n}\{1-o(1)\}.
\]
Taking negative logarithms, we obtain
\[
-\log R_n^*(a,b)
\le
2\lambda m_n+\log 8+o(1).
\]
Since $m_n/L_n\to1+\varepsilon$,
\[
\limsup_{n\to\infty}
\frac{-\log R_n^*(a,b)}{L_n}
\le
2(1+\varepsilon)\lambda.
\]
As $\varepsilon>0$ is arbitrary, we conclude that
\begin{equation}\label{eq:proof-limsup}
\limsup_{n\to\infty}
\frac{-\log R_n^*(a,b)}{L_n}
\le 2\lambda.
\end{equation}

Combining \eqref{eq:proof-liminf}
and \eqref{eq:proof-limsup} establishes
\[
\lim_{n\to\infty}
\frac{-\log R_n^*(a,b)}{L_n}
=
2\lambda
=
2\log\left(
\frac{\sqrt b+\sqrt a}
{\sqrt b-\sqrt a}
\right),
\]
as claimed.
\end{proof}

\begin{lemma}[Uniform covariance estimation]\label{lem:covariance-estimation}
There exist $K_0,K_1\in(0,\infty)$ depending only on $(a,b)$ such that, for all $m,n\ge1$,
\[
\sup_{Q\in\mathcal Q(a,b)}\mathbb E_Q[(\widehat C_{n,m}-C(Q))^2]\le K_0\left(r^{2m}+\frac{K_1^m m!}{n}+\frac1n\right).
\]
\end{lemma}
\begin{proof}
Write $\bar\theta=\overline\theta_{n,m}$, $\widehat\eta=n^{-1}\sum_iX_{i1}$, and $\widehat\mu=n^{-1}\sum_iD_i$. The identity
\[
\widehat C-C=(\widehat\eta-\eta)-\{(\widehat\mu-\mu)\bar\theta+\mu(\bar\theta-\theta)\}
\]
and $0\le\bar\theta\le1$, $\mu\le b$ imply
\[
|\widehat C-C|\le|\widehat\eta-\eta|+|\widehat\mu-\mu|+b|\bar\theta-\theta|.
\]
The variances of $D$ and $X_1$ are bounded by $b+b^2$ by total variance, and projection does not increase squared error relative to $\theta\in[0,1]$. Square the last display, apply $(x+y+z)^2\le3(x^2+y^2+z^2)$, and use Lemma~\ref{lem:minimax-upper}. One admissible choice is $K_1=C_1$ and $K_0=\max\{12b^2,3b^2C_0,6(b+b^2)\}$.
\end{proof}

\begin{proof}[Proof of Theorem~\ref{thm:testing}]
Throughout the proof, let
\[
L_n=\frac{\log n}{\log\log n},
\qquad
r=\frac{\sqrt b-\sqrt a}{\sqrt b+\sqrt a},
\qquad
\lambda=-\log r>0.
\]
All asymptotic statements are understood as
$n\to\infty$, with $0<a<b<\infty$ fixed.

We first establish the existence of consistent
tests above the claimed separation boundary.
Fix $\varepsilon\in(0,1)$ and define
\[
m_n=\lfloor L_n\rfloor,
\qquad
\rho_n=\exp\{-(1-\varepsilon)\lambda L_n\}.
\]
For all sufficiently large $n$, consider
the test
\[
\phi_n=
\mathbf 1\left\{
|\widehat C_{n,m_n}|>\frac{\rho_n}{2}
\right\},
\]
where $\widehat C_{n,m_n}$ is the covariance
estimator defined in Section~\ref{subsec:testing}.

Under the null hypothesis $C(Q)=0$,
the event $\{\phi_n=1\}$ implies
\[
|\widehat C_{n,m_n}-C(Q)|>\rho_n/2.
\]
Under any alternative satisfying
$|C(Q)|\ge\rho_n$, the event
$\{\phi_n=0\}$ implies, by the reverse
triangle inequality,
\[
|\widehat C_{n,m_n}-C(Q)|
\ge \rho_n/2.
\]
Consequently, Markov's inequality yields
\[
\mathcal E_n(\phi_n,\rho_n)
\le
\frac{8}{\rho_n^2}
\sup_{Q\in\mathcal Q(a,b)}
\mathbb E_Q
\left[
(\widehat C_{n,m_n}-C(Q))^2
\right].
\]

By Lemma~\ref{lem:covariance-estimation},
there exist constants $K_0,K_1>0$,
depending only on $(a,b)$, such that
\[
\mathcal E_n(\phi_n,\rho_n)
\le
\frac{8K_0}{\rho_n^2}
\left(
r^{2m_n}
+\frac{K_1^{m_n}m_n!}{n}
+\frac1n
\right).
\]

Since $m_n=L_n+O(1)$, Stirling's formula
gives
\[
\log\left(
\frac{K_1^{m_n}m_n!}{n}
\right)
=
-L_n\log\log\log n+O(L_n).
\]
In particular,
\[
\frac{K_1^{m_n}m_n!}{n}
=o(r^{2m_n}).
\]
Moreover, since
$\log n/L_n=\log\log n\to\infty$,
\[
\frac1n=o(r^{2m_n}).
\]
Therefore,
\[
\begin{aligned}
\mathcal E_n(\phi_n,\rho_n)
&\le
8K_0
\frac{r^{2m_n}}{\rho_n^2}
\{1+o(1)\}\\
&\le
\exp\{-2\varepsilon\lambda L_n+O(1)\}
\longrightarrow0.
\end{aligned}
\]
This proves the first assertion.

We next establish the impossibility
result below the claimed separation
boundary. Fix an arbitrary
$\varepsilon>0$ and set
\[
m_n=
\left\lceil
(1+\varepsilon/2)L_n
\right\rceil,
\qquad
\rho_n=
\exp\{-(1+\varepsilon)\lambda L_n\}.
\]

For each $m\ge1$, consider the signed
measure $\nu_m$ and the probability
distributions $Q_{+,m}$ and $Q_{-,m}$
constructed in the proof of
Lemma~\ref{lem:minimax-lower}.
Recall that
\[
\nu_m=\nu_m^+-\nu_m^-,
\qquad
F_m=\nu_m^++\nu_m^-,
\]
and
\[
\|\nu_m\|_{\mathrm{var}}=1,
\qquad
\int s^j\,\nu_m(ds)=0,
\quad j=1,\ldots,m,
\]
with
\[
\int\nu_m(ds)=E_m
=\frac1{T_m(z)},
\qquad
z=\frac{a+b}{b-a}.
\]

Define
\[
Q_{0,m}
=
\frac12(Q_{+,m}+Q_{-,m}).
\]
By construction,
\[
Q_{0,m}(ds,du)
=
F_m(ds)
\left\{
\frac12\delta_0(du)
+\frac12\delta_1(du)
\right\}.
\]
Thus, under $Q_{0,m}$, the latent
variables $S$ and $U$ are independent,
with $U\sim\operatorname{Bernoulli}(1/2)$.
Consequently,
\[
C(Q_{0,m})
=
\operatorname{Cov}_{Q_{0,m}}(S,U)
=0,
\]
so $Q_{0,m}\in\mathcal Q_0(a,b)$.

We now calculate the covariance under
$Q_{+,m}$. Define
\[
\mu_m=\int_a^b s\,F_m(ds).
\]
Since $F_m$ is supported on $[a,b]$,
\[
a\le\mu_m\le b.
\]
The moment cancellation condition
$\int s\,\nu_m(ds)=0$ implies
\[
\int s\,\nu_m^+(ds)
=
\int s\,\nu_m^-(ds)
=
\frac{\mu_m}{2}.
\]
Therefore,
\[
\mathbb E_{Q_{+,m}}[SU]
=
\int s\,\nu_m^+(ds)
=
\frac{\mu_m}{2}.
\]
Moreover, the Jordan decomposition gives
\[
\mathbb E_{Q_{+,m}}[U]
=
\nu_m^+([a,b])
=
\frac{1+E_m}{2}.
\]
Since $\mathbb E_{Q_{+,m}}[S]=\mu_m$,
we obtain
\[
\begin{aligned}
C(Q_{+,m})
&=
\mathbb E_{Q_{+,m}}[SU]
-\mathbb E_{Q_{+,m}}[S]\,
 \mathbb E_{Q_{+,m}}[U]\\
&=
\frac{\mu_m}{2}
-\mu_m\frac{1+E_m}{2}\\
&=
-\frac{\mu_m E_m}{2}.
\end{aligned}
\]
In particular,
\begin{equation}\label{eq:testing-cov-gap}
|C(Q_{+,m})|
=
\frac{\mu_m E_m}{2}
\ge
\frac a2 E_m.
\end{equation}

Using the Chebyshev identity
\[
E_m=
\frac{2r^m}{1+r^{2m}}
\ge r^m,
\]
we obtain
\[
|C(Q_{+,m})|
\ge \frac a2 r^m.
\]
For $m=m_n$, it follows that
\[
\begin{aligned}
\log\left(
\frac{(a/2)r^{m_n}}{\rho_n}
\right)
&=
\log(a/2)
-\lambda m_n
+(1+\varepsilon)\lambda L_n\\
&=
\frac{\varepsilon\lambda}{2}L_n
+O(1)
\longrightarrow\infty.
\end{aligned}
\]
Consequently, for all sufficiently
large $n$,
\[
|C(Q_{+,m_n})|\ge\rho_n.
\]
Thus $Q_{+,m_n}$ belongs to the
alternative hypothesis, whereas
$Q_{0,m_n}$ belongs to the null
hypothesis.

Let $P_{+,m}$, $P_{-,m}$, and $P_{0,m}$
denote the distributions of one observed
count vector under $Q_{+,m}$,
$Q_{-,m}$, and $Q_{0,m}$, respectively.
The observation law depends linearly
on the mixing distribution, so
\[
P_{0,m}
=
\frac12(P_{+,m}+P_{-,m}).
\]
Hence
\[
d_{\mathrm{TV}}(P_{+,m},P_{0,m})
=
\frac12d_{\mathrm{TV}}(P_{+,m},P_{-,m}).
\]

By the product coupling inequality
and Lemma~\ref{lem:minimax-lower},
\[
\begin{aligned}
&d_{\mathrm{TV}}
\left(
P_{+,m_n}^{\otimes n},
P_{0,m_n}^{\otimes n}
\right)\\
&\quad\le
n\,d_{\mathrm{TV}}
(P_{+,m_n},P_{0,m_n})\\
&\quad=
\frac n2
d_{\mathrm{TV}}
(P_{+,m_n},P_{-,m_n})\\
&\quad\le
\frac{n(2b)^{m_n+1}}
{2(m_n+1)!}.
\end{aligned}
\]

Stirling's formula and
$m_n=(1+\varepsilon/2)L_n+O(1)$
give
\[
\begin{aligned}
&\log\left(
\frac{n(2b)^{m_n+1}}
{2(m_n+1)!}
\right)\\
&\qquad=
-\frac{\varepsilon}{2}\log n
+O(L_n\log\log\log n).
\end{aligned}
\]
Since
$L_n\log\log\log n=o(\log n)$,
we conclude that
\[
d_{\mathrm{TV}}
\left(
P_{+,m_n}^{\otimes n},
P_{0,m_n}^{\otimes n}
\right)
\longrightarrow0.
\]

For any test $\phi_n$, the maximal
sum of testing errors satisfies
\[
\begin{aligned}
\mathcal E_n(\phi_n,\rho_n)
&\ge
P_{0,m_n}^{\otimes n}(\phi_n=1)
+
P_{+,m_n}^{\otimes n}(\phi_n=0)\\
&\ge
1-
d_{\mathrm{TV}}
\left(
P_{+,m_n}^{\otimes n},
P_{0,m_n}^{\otimes n}
\right).
\end{aligned}
\]
Taking the infimum over all tests yields
\[
\liminf_{n\to\infty}
\inf_{\phi_n}
\mathcal E_n(\phi_n,\rho_n)
\ge1.
\]
Conversely, the constant test
$\phi_n\equiv0$ has maximal sum
of errors equal to one for all
sufficiently large $n$, since
the alternative is nonempty.
Therefore,
\[
\limsup_{n\to\infty}
\inf_{\phi_n}
\mathcal E_n(\phi_n,\rho_n)
\le1.
\]
Combining the two inequalities gives
\[
\inf_{\phi_n}
\mathcal E_n\left(
\phi_n,
\exp\{-(1+\varepsilon)\lambda L_n\}
\right)
\longrightarrow1.
\]
This establishes the second assertion
and completes the proof.
\end{proof}

\begin{lemma}[Parametric ratio risk]\label{lem:ratio-risk}
There is $K=K(a,b)<\infty$ such that for all $n\ge1$,
\[
\sup_{Q\in\mathcal Q_0(a,b)}\mathbb E_Q[(\widehat\theta_{R,n}-\theta(Q))^2]\le K/n.
\]
\end{lemma}
\begin{proof}
Put $\mu=\mathbb ED\ge a$, $\eta=\mathbb EX_1$, $\theta_R=\eta/\mu\in[0,1]$. On $\{\widehat\mu\ge a/2\}$,
\[
|\widehat\theta_R-\theta_R|=\frac{|(\widehat\eta-\eta)-\theta_R(\widehat\mu-\mu)|}{\widehat\mu}\le\frac2a(|\widehat\eta-\eta|+|\widehat\mu-\mu|).
\]
Both sample-mean variances are at most $(b+b^2)/n$. On the complementary event, squared error is at most one, and $\Pr(\widehat\mu<a/2)\le4(b+b^2)/(na^2)$. Thus the uniform risk relative to $\theta_R$ is at most $20(b+b^2)/(na^2)$, and $\theta_R=\theta$ on $\mathcal Q_0$.
\end{proof}

\begin{lemma}[Uniform median-of-means concentration]
\label{lem:mom-concentration}
Let $m_n=\max\{1,\lfloor L_n\rfloor\}$,
$V_m=C_0C_1^m m!$, and
$K_n=\lceil32\log n\rceil$.
For all sufficiently large $n$, the partition
into $K_n$ nonempty blocks is possible, and
the projected median of block means obeys
\[
\sup_{Q\in\mathcal Q(a,b)}
\Pr_Q\left(
|\widetilde\theta_{P,n}-\theta(Q)|
>
2r^{m_n}
+17\sqrt{\frac{V_{m_n}\log n}{n}}
\right)
\le n^{-4}.
\]
\end{lemma}

\begin{proof}
For all sufficiently large $n$,
$K_n\le n/2$, and each block contains at
least $n/(2K_n)$ observations.

Since
\[
\operatorname{Var}_Q(Z_{i,m_n})
\le \mathbb E_Q[Z_{i,m_n}^2]
\le V_{m_n},
\]
the variance of each block mean is bounded
uniformly over $Q\in\mathcal Q(a,b)$ by
\[
\frac{2K_nV_{m_n}}{n}.
\]

Set
\[
t_n=\sqrt{\frac{8K_nV_{m_n}}{n}}.
\]
By Chebyshev's inequality, each block mean
deviates from $\mathbb E_Q Z_{i,m_n}$ by
more than $t_n$ with probability at most
$1/4$.

The block means are independent. If their
median differs from $\mathbb E_Q Z_{i,m_n}$
by more than $t_n$, at least half of the
blocks must deviate by more than $t_n$.

Hoeffding's inequality applied to the
independent indicators of these deviations
therefore yields
\[
\Pr_Q\left(
\left|\operatorname{Med}
(\overline Z_1,\ldots,\overline Z_{K_n})
-\mathbb E_Q Z_{i,m_n}\right|>t_n
\right)
\le \exp(-K_n/8)
\le n^{-4}.
\]

Since
\[
K_n=\lceil32\log n\rceil
\le32\log n+1,
\]
for all sufficiently large $n$,
\[
8K_n
\le256\log n+8
\le289\log n.
\]
Consequently,
\[
t_n
=\sqrt{\frac{8K_nV_{m_n}}{n}}
\le17\sqrt{\frac{V_{m_n}\log n}{n}}.
\]

Moreover, the polynomial approximation
bound gives
\[
|\mathbb E_Q Z_{i,m_n}-\theta(Q)|
\le2r^{m_n}.
\]

Combining the preceding inequalities
with the triangle inequality yields
\[
\Pr_Q\left(
\left|
\operatorname{Med}
(\overline Z_1,\ldots,\overline Z_{K_n})
-\theta(Q)
\right|
>
2r^{m_n}
+17\sqrt{\frac{V_{m_n}\log n}{n}}
\right)
\le n^{-4}.
\]

Finally, projection onto $[0,1]$ cannot
increase the distance to
$\theta(Q)\in[0,1]$.
All the bounds hold uniformly over
$Q\in\mathcal Q(a,b)$, which completes
the proof.
\end{proof}

\begin{lemma}[Uniform ratio concentration]
\label{lem:ratio-concentration}
Let
\[
C_b=\frac{eb}{2}+\frac{b^2}{8},
\qquad
A_1=\frac{4\sqrt{24C_b}}{a}.
\]
For all sufficiently large $n$,
\[
\sup_{Q\in\mathcal Q(a,b)}
\Pr_Q\left(
|\widehat\theta_{R,n}-\theta_R(Q)|
>
A_1\sqrt{\frac{\log n}{n}}
\right)
\le n^{-4}.
\]
\end{lemma}

\begin{proof}
Let $W=D$ or $W=X_1$. Conditionally on
its latent Poisson intensity
$\Lambda\in[0,b]$, the variable $W$
has a Poisson distribution with mean
$\Lambda$.

For $|t|\le1$, Hoeffding's lemma
applied to the bounded random variable
$\Lambda$ yields
\[
\begin{aligned}
\log\mathbb E e^{t(W-\mathbb EW)}
&=
\log\mathbb E
\exp\left\{
t(\Lambda-\mathbb E\Lambda)
+\Lambda(e^t-1-t)
\right\}\\
&\le
b(e^t-1-t)+\frac{b^2t^2}{8}\\
&\le C_bt^2,
\end{aligned}
\]
where
\[
C_b=\frac{eb}{2}+\frac{b^2}{8}.
\]

Now let $W_1,\ldots,W_n$ be
independent copies of $W$ and define
\[
\overline W_n=\frac1n\sum_{i=1}^n W_i.
\]

By independence, for every $|t|\le1$,
\[
\log\mathbb E
\exp\left\{
nt(\overline W_n-\mathbb EW)
\right\}
\le nC_bt^2.
\]

Applying Chernoff's inequality and
choosing
\[
t=\frac{x}{2C_b}
\]
gives, for $0<x\le2C_b$,
\[
\Pr\left(
\overline W_n-\mathbb EW>x
\right)
\le
\exp\left(
-\frac{nx^2}{4C_b}
\right).
\]

The same argument applied to
$-\overline W_n$ yields
\[
\Pr\left(
|\overline W_n-\mathbb EW|>x
\right)
\le
2\exp\left(
-\frac{nx^2}{4C_b}
\right).
\]

Set
\[
x_n=\sqrt{\frac{24C_b\log n}{n}}.
\]
For all sufficiently large $n$,
\[
x_n\le\min\{2C_b,a/2\}.
\]

Applying the preceding concentration
inequality to $W=X_1$ and $W=D$,
and using the union bound, gives
\[
\Pr_Q\left(
|\widehat\eta-\eta|>x_n
\ \text{or}\
|\widehat\mu-\mu|>x_n
\right)
\le4n^{-6}.
\]

On the complementary event,
\[
\widehat\mu
\ge\mu-x_n
\ge a/2.
\]

Since $\theta_R=\eta/\mu\in[0,1]$,
we obtain
\[
\begin{aligned}
|\widehat\theta_{R,n}-\theta_R|
&=
\frac{
|(\widehat\eta-\eta)
-\theta_R(\widehat\mu-\mu)|
}{\widehat\mu}\\
&\le
\frac2a
\left(
|\widehat\eta-\eta|
+
|\widehat\mu-\mu|
\right)\\
&\le\frac{4x_n}{a}\\
&=
A_1\sqrt{\frac{\log n}{n}}.
\end{aligned}
\]

Consequently,
\[
\sup_{Q\in\mathcal Q(a,b)}
\Pr_Q\left(
|\widehat\theta_{R,n}-\theta_R(Q)|
>
A_1\sqrt{\frac{\log n}{n}}
\right)
\le4n^{-6}
\le n^{-4}
\]
for all sufficiently large $n$.
\end{proof}

\begin{proof}[Proof of Theorem~\ref{thm:adaptive}]
Let
\[
L_n=\frac{\log n}{\log\log n},
\qquad
m_n=\max\{1,\lfloor L_n\rfloor\},
\qquad
\lambda=-\log r>0.
\]
Throughout the proof, the support bounds
$0<a<b<\infty$ are fixed.

Let
\[
V_m=C_0C_1^m m!,
\qquad
c=32,\qquad A_0=17,
\qquad
A_1=\frac{4\sqrt{24C_b}}{a},
\]
where $C_0,C_1>0$ are the constants
from Lemma~\ref{lem:minimax-upper} and
\[
C_b=\frac{eb}{2}+\frac{b^2}{8}.
\]
All these constants depend only on the
fixed support bounds $(a,b)$.

For sufficiently large $n$, define
\[
K_n=\lceil32\log n\rceil.
\]
We restrict attention to sample sizes
large enough that
\[
K_n\le\frac n2.
\]
Thus the sample can be partitioned into
$K_n$ nonempty blocks whose sizes differ
by at most one, as required for the
median-of-means construction of
$\widetilde\theta_{P,n}$.

For the finitely many smaller sample
sizes, define
$\widehat\theta_n^{\mathrm{ad}}
=\widehat\theta_{R,n}$.
This convention makes the estimator
well-defined for every $n$ and has no
effect on the asymptotic risk statements.

For sufficiently large $n$, set
\[
e_n=
2r^{m_n}
+17\sqrt{\frac{V_{m_n}\log n}{n}},
\qquad
h_n=
A_1\sqrt{\frac{\log n}{n}}.
\]
For the remaining sample sizes,
$e_n$ and $h_n$ may be assigned
arbitrary positive deterministic values.
The adaptive estimator is then given
by the comparison rule
\eqref{eq:adaptive-selector} for
sufficiently large $n$.

For a generic $Q\in\mathcal Q(a,b)$,
write
\[
\theta=\theta(Q),
\qquad
\theta_R=
\frac{\mathbb E_Q[SU]}{\mathbb E_Q[S]}.
\]
Consider the events
\[
G_P=
\left\{
|\widetilde\theta_{P,n}-\theta|
\le e_n
\right\},
\qquad
G_R=
\left\{
|\widehat\theta_{R,n}-\theta_R|
\le h_n
\right\}.
\]

By Lemmas~\ref{lem:mom-concentration}
and~\ref{lem:ratio-concentration},
for all sufficiently large $n$,
\begin{equation}\label{eq:adaptive-good-events}
\begin{aligned}
\sup_{Q\in\mathcal Q(a,b)}
P_Q(G_P^c)&\le n^{-4},\\
\sup_{Q\in\mathcal Q(a,b)}
P_Q(G_R^c)&\le n^{-4}.
\end{aligned}
\end{equation}
The constants in these concentration
inequalities depend only on the fixed
support bounds $(a,b)$.
In particular, the bounds hold uniformly
over all $Q\in\mathcal Q(a,b)$,
without any independence or smoothness
assumption on the latent variables.

We first establish the parametric risk
bound on the zero-covariance submodel.

Let $Q\in\mathcal Q_0(a,b)$.
Since $\operatorname{Cov}_Q(S,U)=0$,
the identity
\[
\theta_R-\theta
=
\frac{\operatorname{Cov}_Q(S,U)}
{\mathbb E_Q[S]}
\]
implies $\theta_R=\theta$.

On $G_P\cap G_R$, the triangle
inequality gives
\[
\begin{aligned}
|\widetilde\theta_{P,n}
-\widehat\theta_{R,n}|
&\le
|\widetilde\theta_{P,n}-\theta|
+
|\widehat\theta_{R,n}-\theta|\\
&\le e_n+h_n.
\end{aligned}
\]
Thus, by the selection rule,
\[
\widehat\theta_n^{\mathrm{ad}}
=\widehat\theta_{R,n}
\qquad\text{on }G_P\cap G_R.
\]

Since the two estimators and the
target take values in $[0,1]$,
their squared errors are bounded
by one. Therefore,
\[
\begin{aligned}
&\mathbb E_Q\left[
(\widehat\theta_n^{\mathrm{ad}}
-\theta)^2\right]\\
&\quad\le
\mathbb E_Q\left[
(\widehat\theta_{R,n}
-\theta)^2\right]
+
P_Q(G_P^c\cup G_R^c).
\end{aligned}
\]
Using \eqref{eq:adaptive-good-events}
and the uniform ratio risk bound
from Lemma~\ref{lem:ratio-risk},
we obtain
\[
\begin{aligned}
&\sup_{Q\in\mathcal Q_0(a,b)}
\mathbb E_Q\left[
(\widehat\theta_n^{\mathrm{ad}}
-\theta(Q))^2\right]\\
&\quad\le
\frac{K_{a,b}}{n}+2n^{-4}
=O(n^{-1}),
\end{aligned}
\]
where $K_{a,b}>0$ depends only
on $(a,b)$.
This proves the first assertion.

We next establish the logarithmic
risk bound over the unrestricted
class $\mathcal Q(a,b)$.

Fix an arbitrary
$Q\in\mathcal Q(a,b)$.
On $G_P$, if the polynomial
estimator is selected, then
\[
|\widehat\theta_n^{\mathrm{ad}}
-\theta|
=
|\widetilde\theta_{P,n}-\theta|
\le e_n.
\]

Otherwise, the ratio estimator is
selected, and the selection rule
implies
\[
|\widehat\theta_{R,n}
-\widetilde\theta_{P,n}|
\le e_n+h_n.
\]
Consequently, on $G_P$,
\[
\begin{aligned}
|\widehat\theta_n^{\mathrm{ad}}
-\theta|
&=
|\widehat\theta_{R,n}-\theta|\\
&\le
|\widehat\theta_{R,n}
-\widetilde\theta_{P,n}|
+
|\widetilde\theta_{P,n}-\theta|\\
&\le2e_n+h_n.
\end{aligned}
\]
Thus, regardless of which estimator
is selected,
\[
|\widehat\theta_n^{\mathrm{ad}}-\theta|
\le2e_n+h_n
\qquad\text{on }G_P.
\]

On $G_P^c$, the squared error is
at most one. It follows that
\[
\begin{aligned}
&\mathbb E_Q\left[
(\widehat\theta_n^{\mathrm{ad}}
-\theta)^2\right]\\
&\quad\le
(2e_n+h_n)^2+P_Q(G_P^c).
\end{aligned}
\]
Taking the supremum over
$\mathcal Q(a,b)$ gives
\begin{equation}\label{eq:adaptive-global-bound}
\begin{aligned}
&\sup_{Q\in\mathcal Q(a,b)}
\mathbb E_Q\left[
(\widehat\theta_n^{\mathrm{ad}}
-\theta(Q))^2\right]\\
&\quad\le
(2e_n+h_n)^2+n^{-4}.
\end{aligned}
\end{equation}

It remains to determine the
asymptotic order of the right-hand
side. For all sufficiently large
$n$, $m_n=\lfloor L_n\rfloor$,
so that $m_n=L_n+O(1)$.

By Stirling's formula,
\[
\begin{aligned}
\log\left(
\frac{V_{m_n}\log n}{n}
\right)
&=
\log C_0+m_n\log C_1
+\log(m_n!)\\
&\quad+\log\log n-\log n\\
&=
-L_n\log\log\log n+O(L_n).
\end{aligned}
\]
Since
\[
\log(r^{2m_n})
=-2\lambda L_n+O(1),
\]
we obtain
\[
\frac{V_{m_n}\log n}{n}
=o(r^{2m_n}).
\]
Therefore,
\[
e_n=2r^{m_n}+o(r^{m_n}).
\]

Moreover,
\[
\log h_n
=-\frac12\log n+O(\log\log n),
\]
whereas
\[
\log(r^{m_n})
=-\lambda L_n+O(1).
\]
Since
$\log n/L_n=\log\log n\to\infty$,
it follows that
\[
h_n=o(r^{m_n}).
\]

Combining these estimates yields
\[
2e_n+h_n
=
4r^{m_n}+o(r^{m_n}).
\]
Substituting into
\eqref{eq:adaptive-global-bound},
and observing that
$n^{-4}=o(r^{2m_n})$, gives
\[
\begin{aligned}
&\sup_{Q\in\mathcal Q(a,b)}
\mathbb E_Q\left[
(\widehat\theta_n^{\mathrm{ad}}
-\theta(Q))^2\right]\\
&\qquad\le
(16+o(1))r^{2m_n}.
\end{aligned}
\]
Consequently,
\[
\liminf_{n\to\infty}
\frac{
-\log\displaystyle\sup_{Q\in\mathcal Q(a,b)}
\mathbb E_Q\left[
(\widehat\theta_n^{\mathrm{ad}}
-\theta(Q))^2\right]
}{L_n}
\ge2\lambda.
\]

This inequality follows solely from
the risk analysis of the adaptive
estimator. To establish equality
of the logarithmic exponent, we
invoke the minimax lower bound.

Indeed, by the definition of the
minimax risk,
\[
R_n^*(a,b)
\le
\sup_{Q\in\mathcal Q(a,b)}
\mathbb E_Q\left[
(\widehat\theta_n^{\mathrm{ad}}
-\theta(Q))^2\right].
\]
The sharp minimax result from
Theorem~\ref{thm:minimax-rate}
therefore implies
\[
\limsup_{n\to\infty}
\frac{
-\log\displaystyle\sup_{Q\in\mathcal Q(a,b)}
\mathbb E_Q\left[
(\widehat\theta_n^{\mathrm{ad}}
-\theta(Q))^2\right]
}{L_n}
\le2\lambda.
\]
Combining the two inequalities yields
\[
\lim_{n\to\infty}
\frac{
-\log\displaystyle\sup_{Q\in\mathcal Q(a,b)}
\mathbb E_Q\left[
(\widehat\theta_n^{\mathrm{ad}}
-\theta(Q))^2\right]
}{L_n}
=2\lambda.
\]

Thus the adaptive construction achieves
parametric squared-error risk on the
zero-covariance submodel and, together
with the minimax lower bound, attains
the optimal full-model logarithmic
exponent.
\end{proof}


\begin{thebibliography}{99}
\bibitem[Aitchison and Ho(1989)]{AitchisonHo1989}
Aitchison, J. and Ho, C. H. (1989).
The multivariate Poisson-log normal distribution.
\emph{Biometrika}, \textbf{76}(4), 643--653.
\url{https://doi.org/10.1093/biomet/76.4.643}.

\bibitem[Baraud(2002)]{Baraud2002}
Baraud, Y. (2002).
Non-asymptotic minimax rates of testing in signal detection.
\emph{Bernoulli}, \textbf{8}(5), 577--606.

\bibitem[Boucheron et al.(2013)Boucheron, Lugosi, and Massart]{BoucheronLugosiMassart2013}
Boucheron, S., Lugosi, G. and Massart, P. (2013).
\emph{Concentration Inequalities: A Nonasymptotic Theory of Independence}.
Oxford University Press.

\bibitem[Cai and Low(2004)]{CaiLow2004}
Cai, T. T. and Low, M. G. (2004).
Minimax estimation of linear functionals over nonconvex parameter spaces.
\emph{The Annals of Statistics}, \textbf{32}(2), 552--576.
\url{https://doi.org/10.1214/009053604000000094}.

\bibitem[Carroll and Hall(1988)]{CarrollHall1988}
Carroll, R. J. and Hall, P. (1988).
Optimal rates of convergence for deconvolving a density.
\emph{Journal of the American Statistical Association},
\textbf{83}(404), 1184--1186.

\bibitem[Catoni(2012)]{Catoni2012}
Catoni, O. (2012).
Challenging the empirical mean and empirical variance: a deviation study.
\emph{Annales de l'Institut Henri Poincar\'e, Probabilit\'es et Statistiques},
\textbf{48}(4), 1148--1185.
\url{https://doi.org/10.1214/11-AIHP454}.

\bibitem[Cavalier and Hengartner(2009)]{CavalierHengartner2009}
Cavalier, L. and Hengartner, N. W. (2009).
Estimating Linear Functionals in Poisson Mixture Models.
\emph{Journal of Nonparametric Statistics}, \textbf{21}(6), 713--728.
\url{https://doi.org/10.1080/10485250902971716}.

\bibitem[Cheney(1966)]{Cheney1966}
Cheney, E. W. (1966).
\emph{Introduction to Approximation Theory}.
McGraw-Hill, New York.

\bibitem[DeVore and Lorentz(1993)]{DeVoreLorentz1993}
DeVore, R. A. and Lorentz, G. G. (1993).
\emph{Constructive Approximation}.
Grundlehren der mathematischen Wissenschaften, Vol.~303.
Springer.

\bibitem[Devroye et al.(2016)Devroye, Lerasle, Lugosi, and Oliveira]{DevroyeLerasleLugosiOliveira2016}
Devroye, L., Lerasle, M., Lugosi, G. and Oliveira, R. I. (2016).
Sub-Gaussian Mean Estimators.
\emph{The Annals of Statistics}, \textbf{44}(6), 2695--2725.
\url{https://doi.org/10.1214/16-AOS1440}.

\bibitem[Donoho(1994)]{Donoho1994}
Donoho, D. L. (1994).
Statistical Estimation and Optimal Recovery.
\emph{The Annals of Statistics}, \textbf{22}(1), 238--270.
\url{https://doi.org/10.1214/aos/1176325367}.

\bibitem[Fan(1991)]{Fan1991}
Fan, J. (1991).
On the optimal rates of convergence for nonparametric deconvolution problems.
\emph{The Annals of Statistics}, \textbf{19}(3), 1257--1272.

\bibitem[Favaro et al.(2026)Favaro, Fortini, and Jana]{FavaroFortiniJana2026}
Favaro, S., Fortini, S. and Jana, S. (2026).
Poisson Empirical Bayes Estimation of Sums of Random Variables
via Minimum-Distance Methods.
\emph{arXiv preprint} arXiv:2610.07190 [stat.ME].
\url{https://arxiv.org/abs/2610.07190}.

\bibitem[Feller(1968)]{Feller1968}
Feller, W. (1968).
\emph{An Introduction to Probability Theory and Its Applications},
Vol.~1, 3rd ed.
John Wiley \& Sons, New York.
ISBN 9780471257080.

\bibitem[Goldenshluger and Lepski(2011)]{GoldenshlugerLepski2011}
Goldenshluger, A. and Lepski, O. (2011).
Bandwidth selection in kernel density estimation:
oracle inequalities and adaptive minimax optimality.
\emph{The Annals of Statistics}, \textbf{39}(3), 1608--1632.
\url{https://doi.org/10.1214/11-AOS883}.

\bibitem[Hengartner(1997)]{Hengartner1997}
Hengartner, N. W. (1997).
Adaptive Demixing in Poisson Mixture Models.
\emph{The Annals of Statistics}, \textbf{25}(3), 917--928.
\url{https://doi.org/10.1214/aos/1069362730}.

\bibitem[Holgate(1964)]{Holgate1964}
Holgate, P. (1964).
Estimation for the Bivariate Poisson Distribution.
\emph{Biometrika}, \textbf{51}(1/2), 241--245.
\url{https://doi.org/10.2307/2334210}.

\bibitem[Ingster and Suslina(2003)]{IngsterSuslina2003}
Ingster, Y. I. and Suslina, I. A. (2003).
\emph{Nonparametric Goodness-of-Fit Testing Under Gaussian Models}.
Lecture Notes in Statistics, Vol.~169.
Springer, New York.
ISBN 9780387955315.
\url{https://doi.org/10.1007/978-0-387-21580-8}.

\bibitem[Jana et al.(2025)Jana, Polyanskiy, and Wu]{OptimalEB2025}
Jana, S., Polyanskiy, Y. and Wu, Y. (2025).
Optimal Empirical Bayes Estimation for the Poisson Model
via Minimum-Distance Methods.
\emph{Information and Inference: A Journal of the IMA},
\textbf{14}(4), iaaf027.
\url{https://doi.org/10.1093/imaiai/iaaf027}.

\bibitem[Karlis and Xekalaki(2005)]{KarlisXekalaki2005}
Karlis, D. and Xekalaki, E. (2005).
Mixed Poisson Distributions.
\emph{International Statistical Review}, \textbf{73}(1), 35--58.
\url{https://doi.org/10.1111/j.1751-5823.2005.tb00250.x}.

\bibitem[Kiefer and Wolfowitz(1956)]{KieferWolfowitz1956}
Kiefer, J. and Wolfowitz, J. (1956).
Consistency of the maximum likelihood estimator in the presence
of infinitely many incidental parameters.
\emph{The Annals of Mathematical Statistics}, \textbf{27}(4), 887--906.
\url{https://doi.org/10.1214/aoms/1177728066}.

\bibitem[Kingman(1993)]{Kingman1993}
Kingman, J. F. C. (1993).
\emph{Poisson Processes}.
Oxford Studies in Probability, Vol.~3.
Oxford University Press, Oxford.
ISBN 9780198536932.

\bibitem[Laurent and Massart(2000)]{LaurentMassart2000}
Laurent, B. and Massart, P. (2000).
Adaptive estimation of a quadratic functional by model selection.
\emph{The Annals of Statistics}, \textbf{28}(5), 1302--1338.

\bibitem[Lepskii(1991)]{Lepski1991}
Lepskii, O. V. (1991).
On a Problem of Adaptive Estimation in Gaussian White Noise.
\emph{Theory of Probability \& Its Applications}, \textbf{35}(3), 454--466.
\url{https://doi.org/10.1137/1135065}.

\bibitem[Lim and Han(2024)]{LimHan2024}
Lim, K. and Han, F. (2024).
Smoothed NPMLEs in Nonparametric Poisson Mixtures and Beyond.
\emph{arXiv preprint} arXiv:2406.08808 [math.ST].
\url{https://arxiv.org/abs/2406.08808}.

\bibitem[Lindsay(1983)]{Lindsay1983}
Lindsay, B. G. (1983).
The Geometry of Mixture Likelihoods: A General Theory.
\emph{The Annals of Statistics}, \textbf{11}(1), 86--94.
\url{https://doi.org/10.1214/aos/1176346059}.

\bibitem[Lindsay(1995)]{Lindsay1995}
Lindsay, B. G. (1995).
\emph{Mixture Models: Theory, Geometry and Applications}.
NSF-CBMS Regional Conference Series in Probability and Statistics, Vol.~5.
Institute of Mathematical Statistics, Hayward, CA.

\bibitem[Lugosi and Mendelson(2019)]{LugosiMendelson2019}
Lugosi, G. and Mendelson, S. (2019).
Sub-Gaussian estimators of the mean of a random vector.
\emph{The Annals of Statistics}, \textbf{47}(2), 783--794.
\url{https://doi.org/10.1214/17-AOS1639}.

\bibitem[Teicher(1960)]{Teicher1960}
Teicher, H. (1960).
On the Mixture of Distributions.
\emph{The Annals of Mathematical Statistics}, \textbf{31}(1), 55--73.
\url{https://doi.org/10.1214/aoms/1177705987}.

\bibitem[Teicher(1961)]{Teicher1961}
Teicher, H. (1961).
Identifiability of Mixtures.
\emph{The Annals of Mathematical Statistics}, \textbf{32}(1), 244--248.
\url{https://doi.org/10.1214/aoms/1177705155}.

\bibitem[Tsybakov(2009)]{Tsybakov2009}
Tsybakov, A. B. (2009).
\emph{Introduction to Nonparametric Estimation}.
Springer Series in Statistics.
Springer, New York.
\url{https://doi.org/10.1007/b13794}.

\end{thebibliography}
\end{document}